\documentclass[a4paper,oneside,12pt,reqno]{amsart}
\usepackage[utf8]{inputenc}
\usepackage{stmaryrd}
\usepackage{bbm}
\usepackage{xcolor}
\usepackage{mathrsfs}
\usepackage{enumerate}
\usepackage{centernot}
\usepackage{latexsym,amsxtra}
\usepackage{dsfont}
\usepackage{slashed}
\usepackage[all]{xy}
\usepackage{amscd,graphics}
\usepackage{amsmath,amsfonts,amsthm,amssymb}
\usepackage{latexsym,amsmath}
\usepackage{graphicx,psfrag}
\usepackage{hyperref}
\usepackage{fancyhdr}

\usepackage{libertine}
\usepackage{newtxmath}
\fancypagestyle{plain}{%
  \fancyhf{}%
  \fancyhead[R]{\thepage}%
}
\fancypagestyle{firstpage}{%
  \fancyhf{}%
  \fancyhead[R]{\thepage}%
}
\usepackage{cleveref}
\usepackage{comment} 
\usepackage{braket}
\usepackage{caption}
\usepackage{enumitem}
\theoremstyle{plain}
\newtheorem{thm}{Theorem}[section]
\crefname{thm}{Theorem}{Theorems}
\Crefname{thm}{Theorem}{Theorems}
\newtheorem{prop}[thm]{Proposition}
\crefname{prop}{Proposition}{Propositions}
\Crefname{prop}{Proposition}{Propositions}
\newtheorem{lem}[thm]{Lemma}
\crefname{lem}{Lemma}{Lemmas}
\Crefname{lem}{Lemma}{Lemmas}
\newtheorem{cor}[thm]{Corollary}
\crefname{cor}{Corollary}{Corollaries}
\Crefname{cor}{Corollary}{Corollaries}

\crefname{conj}{Conjecture}{Conjectures}
\Crefname{conj}{Conjecture}{Conjectures}

\crefname{cons}{Construction}{Constructions}
\Crefname{cons}{Construction}{Constructions}

\crefname{claim}{Claim}{Claims}
\Crefname{claim}{Claim}{Claims}

\crefname{property}{Property}{Properties}
\Crefname{property}{Property}{Properties}

\crefname{problem}{Problem}{Problems}
\Crefname{problem}{Problem}{Problems}

\theoremstyle{definition}
\newtheorem{defi}[thm]{Definition}
\crefname{defi}{Definition}{Definitions}
\Crefname{defi}{Definition}{Definitions}

\crefname{nota}{Notation}{Notations}
\Crefname{nota}{Notation}{Notations}

\crefname{convention}{Convention}{Conventions}
\Crefname{convention}{Convention}{Conventions}

\crefname{cond}{Condition}{Conditions}
\Crefname{cond}{Condition}{Conditions}

\crefname{assum}{Assumption}{Assumptions}
\Crefname{assum}{Assumption}{Assumptions}

\theoremstyle{remark}
\newtheorem{rem}[thm]{Remark}
\crefname{rem}{Remark}{Remarks}
\Crefname{rem}{Remark}{Remarks}
\newtheorem{ex}[thm]{Example}
\crefname{ex}{Example}{Examples}
\Crefname{ex}{Example}{Examples}

\newtheorem{ques}[thm]{Question}
\crefname{ques}{Question}{Questions}
\Crefname{ques}{Question}{Questions}

\crefname{section}{Section}{Sections}
\Crefname{section}{Section}{Sections}
\crefname{subsection}{Subsection}{Subsections}
\Crefname{subsection}{Subsection}{Subsections}
\crefname{figure}{Figure}{Figures}
\Crefname{figure}{Figure}{Figures}

\newcommand{\Diff}{\textnormal{Diff}}
\newcommand{\BDiff}{B\textnormal{Diff}}
\newcommand{\TDiff}{\textnormal{TDiff}}

\newcommand{\Homeo}{\textnormal{Homeo}}

\newcommand{\ab}{\textnormal{ab}}

\newcommand{\SW}{\mathrm{SW}}
\newcommand{\BF}{\mathrm{BF}}
\newcommand{\FBF}{\mathrm{FBF}}

\newcommand{\id}{\textnormal{id}}

\newcommand{\Z}{\mathbb{Z}}

\newcommand{\R}{\mathbb{R}}
\newcommand{\C}{\mathbb{C}}

\newcommand{\fraks}{\mathfrak{s}}
\newcommand{\frakt}{\mathfrak{t}}

\newcommand{\ind}{\mathrm{ind}}

\DeclareMathOperator{\Map}{Map}

\newcommand{\del}{\partial}

\newcommand{\divi}{\mathrm{div}}

\newcommand{\Dsla}{\slashed{D}}

\newcommand{\pt}{\mathrm{pt}}

\newcommand{\even}{\mathrm{even}}
\newcommand{\odd}{\mathrm{odd}}

\newcommand{\std}{\mathrm{std}}

\newcommand{\Ker}{\mathop{\mathrm{Ker}}\nolimits}
\newcommand{\Coker}{\mathop{\mathrm{Coker}}\nolimits}

\newcommand{\rank}{\mathop{\mathrm{rank}}\nolimits}

\newcommand{\Th}{\operatorname{Th}}
\newcommand{\Pin}{\operatorname{Pin}}
\newcommand{\SO}{\operatorname{SO}}
\newcommand{\GL}{\operatorname{GL}}

\newcommand{\Spin}{\operatorname{Spin}}

\title{Exotic Dehn twists on sums of 4-manifolds, nuclei, and exotic structures}

\author{Hokuto Konno}
\address{Graduate School of Mathematical Sciences, the University of Tokyo, 3-8-1 Komaba, Meguro, Tokyo 153-8914, Japan \\and\\
RIKEN iTHEMS, Wako, Saitama 351-0198, Japan
\\and\\
Theoretical Sciences Visiting Program, Okinawa Institute of Science and Technology Graduate University, Onna, 904-0495, Japan}
\email{konno@ms.u-tokyo.ac.jp}

\dedicatory{In memory of Mikio Furuta}

\begin{document}

\begin{abstract}
For closed 4-manifolds $X_1, X_2$, we give a sufficient condition, in terms of the Seiberg--Witten invariants of $X_i$, for the Dehn twist on $X_1\#X_2$ along the neck to be smoothly nontrivial.
This provides a broad generalization of the result of Kronheimer--Mrowka for $K3\#K3$.
We also show that the boundary Dehn twists on the punctured connected sums of many pairs of elliptic surfaces are smoothly nontrivial.
We also detect new examples of exotic Dehn twists along Seifert fibered 3-manifolds in elliptic surfaces.
As a consequence, we exhibit a new phenomenon in which ``exotic diffeomorphisms detect exotic smooth structures on 4-manifolds": the triviality of the boundary Dehn twist characterizes the standard nucleus $N(n)$ among the nuclei $N(n)_{p,q}$ homeomorphic to it, provided that $n\equiv0 \bmod 4$.
These results are proved by computing the families Bauer--Furuta invariant for a new framing of the families Dirac index.
\end{abstract}

\maketitle
\thispagestyle{empty}

%\tableofcontents

\section{Main results}

\subsection{Summary}
A diffeomorphism $f : X \to X$ of a smooth manifold $X$ is called an \emph{exotic diffeomorphism} if $f$ is topologically isotopic to the identity but not smoothly isotopic to the identity.
Since Ruberman's construction of the first examples of exotic diffeomorphisms of 4-manifolds \cite{Rub98}, the detection and study of such diffeomorphisms have attracted growing attention in four-dimensional topology, particularly in recent years.
Much of the recent work on exotic diffeomorphisms has focused on Dehn twists along embedded 3-spheres or Seifert fibered 3-manifolds in 4-manifolds.
(See Subsections~\ref{subsection Boundary Dehn twists of punctured 4-manifolds}, \ref{subsection Seifert fibered Dehn twists on elliptic surfaces} for the definition of Dehn twists.)

In this paper, by computing the families Bauer--Furuta invariant for a new framing of the families Dirac index, we establish the following results on (i) Dehn twists along 3-spheres, and (ii) Dehn twists along Seifert fibered 3-manifolds.
The results in (ii) provide a new way to detect exotic smooth structures on 4-manifolds through exotic diffeomorphisms.

\begin{itemize}
\item[(i)] For closed 4-manifolds $X_1, X_2$, we give a sufficient condition for the Dehn twist on $X_1\#X_2$ along the neck to be smoothly nontrivial (Theorem~\ref{thm: eg Dehn twists on connected sums}).
This provides a vast generalization of the pioneering example of Kronheimer--Mrowka \cite{KM20Dehn} for $K3\#K3$.
We also give a sufficient condition for the boundary Dehn twist on a punctured 4-manifold to be smoothly nontrivial (Theorem~\ref{thm: eg puncture boundary Dehn twist}).
Connected sums of elliptic surfaces provide a rich class of examples to which these theorems apply.
%We also show that the boundary Dehn twists on the punctured connected sums of many pairs of elliptic surfaces are smoothly nontrivial (Theorem~\ref{thm: eg puncture boundary Dehn twist}).

\item[(ii)] We detect new examples of exotic Dehn twists along Seifert fibered 3-manifolds in elliptic surfaces (Theorem~\ref{thm: Dehn twists on elliptic surfaces}), extending previous results on the detection of Dehn twists on irreducible 4-manifolds~\cite{konno2025constraintslefschetzfibrationsfourdimensional}.
As a consequence, we exhibit a new phenomenon in which ``exotic diffeomorphisms detect exotic smooth structures on 4-manifolds": the triviality of the boundary Dehn twist characterizes the standard nucleus $N(n)$ among the nuclei $N(n)_{p,q}$ homeomorphic to it, provided that $n\equiv0 \bmod 4$ (Corollary~\ref{cor: nuclei}).
\end{itemize}

In the following subsections, we describe these results in more detail.

\subsection{Dehn twists on connected sums}

The first example of an exotic diffeomorphism arising from a Dehn twist was constructed by Kronheimer and Mrowka~\cite{KM20Dehn}, who proved that the Dehn twist on $K3\#K3$ along the $S^3$ in the separating neck is exotic.
However, relatively few examples of exotic Dehn twists along $S^3$ in closed 4-manifolds are known.
To the best of our knowledge, these are limited to the Dehn twist on $K3\#K3$, its extension to $K3\#K3\#(S^2\times S^2)$ by Lin~\cite{Lin20Dehn}, and variants obtained by replacing $K3$ surfaces with homotopy $K3$ surfaces.
Our first theorem provides a fairly general sufficient condition for the Dehn twist along the neck of a connected sum $X_1\#X_2$ to be smoothly nontrivial:

\begin{thm}
\label{thm: eg Dehn twists on connected sums}
Let $X_1,X_2$ be closed simply-connected smooth spin 4-manifolds with $\sigma(X_i)\equiv 16 \bmod 32$. 
For $i=1,2$, let $\fraks_i$ be a spin$^c$ structure on $X_i$ such that $c_1(\fraks_i)$ is divisible by 8.
Assume that $d(\fraks_i)=0$ and  $\SW(X_i,\fraks_i)$ is odd.
Then the Dehn twist $\delta$ on $X_1\#X_2$ along the neck is an exotic diffeomorphism.
\end{thm}

This theorem substantially generalizes the result of Kronheimer--Mrowka~\cite{KM20Dehn} for $K3\#K3$.
To illustrate the scope of this theorem, we apply it to connected sums of elliptic surfaces.
We use the standard notation $E(n)_{p,q}$ for the simply-connected elliptic surface of degree $n$ obtained by logarithmic transformations of multiplicities $p$ and $q$.
For details, see Section~\ref{section Examples}.

\begin{ex}[See Theorem~\ref{thm: elliptic surface sum}]
Let $k,l\geq0$, and let $p,q,r,s$ be positive odd integers such that
$\gcd(p,q)=\gcd(r,s)=1$.
Suppose that
\[
\left\{
\begin{aligned}
&\text{either } k \text{ is even or }
\max\{p,q\}\geq5, \text{and}\\
&\text{either } l \text{ is even or }
\max\{r,s\}\geq5.
\end{aligned}
\right.
\]
Then the Dehn twist $\delta$ on
\[
E(4k+2)_{p,q}\#E(4l+2)_{r,s}
\]
along the neck is an exotic diffeomorphism.
In particular, since $E(n)_{1,1}=E(n)$, the Dehn twist along the neck of 
\[
E(8k+2)\#E(8l+2)
\]
is exotic for all integers $k,l\geq0$.
The case $k=l=0$ recovers the example of Kronheimer--Mrowka~\cite{KM20Dehn} on $E(2)\#E(2)=K3\#K3$.
\end{ex}

\subsection{Boundary Dehn twists on punctured 4-manifolds}

Let $W$ be a smooth 4-manifold with boundary. 
Let $\Diff_\del(W)$ denote the group of diffeomorphisms of $W$ that are the identity near the boundary, and let $\Homeo_\del(W)$ denote the corresponding group of homeomorphisms.
A diffeomorphism $f\in\Diff_\del(W)$ is called a \emph{relatively exotic diffeomorphism} if it is topologically isotopic to the identity through $\Homeo_\del(W)$ but not smoothly isotopic to the identity through $\Diff_\del(W)$.
Given a closed 4-manifold $X$, let $\mathring{X}$ denote the punctured 4-manifold obtained from $X$, namely, $\mathring{X}=X \setminus \mathrm{Int}(D^4)$.
On $\mathring{X}$, we consider the \emph{boundary Dehn twist} $\delta\in\Diff_\del(\mathring{X})$, which is supported in a collar neighborhood of $\del\mathring{X}$ and is given there by the standard model Dehn twist on $S^3\times[0,1]$.

Theorem~\ref{thm: eg Dehn twists on connected sums} immediately implies that the boundary Dehn twist on $\mathring{X}_1$ is nontrivial in $\pi_0(\Diff_\del(\mathring{X}_1))$ for $X_1$ in the statement.
For example, the boundary Dehn twist $\delta$ on the punctured elliptic surface $(E(8k+2))^\circ$ is relatively exotic for all integers $k\geq0$ (see Example~\ref{ex: boundary Dehn E(8k+2)}).

%The above example also strengthens a result of Baraglia and the author~\cite[Theorem~1.5]{baraglia2024irreducible4manifoldsadmitexotic}, which established the nontriviality of boundary Dehn twists on $(E(4k+2)_{p,q})^\circ$ for all but finitely many pairs of relatively prime positive odd integers $(p,q)$.
%In the subfamily $E(8k+2)$, \cite[Theorem~1.4]{baraglia2024irreducible4manifoldsadmitexotic} applies when $k\equiv0\bmod4$. Example~\ref{ex: boundary Dehn E(8k+2)} covers the remaining cases $k\not\equiv0\bmod4$.

Another class of examples, typically arising from connected sums of irreducible 4-manifolds, can be detected by the following sufficient condition in terms of the Bauer--Furuta invariant:

\begin{thm}
\label{thm: eg puncture boundary Dehn twist}
Let $X$ be a simply-connected closed smooth spin 4-manifold with $\sigma(X)\equiv 16$ mod 32.
Suppose that there exists a spin$^c$ structure $\fraks$ on $X$ for which $c_1(\fraks)$ is divisible by 8.
Assume that $d(\fraks)=1$ and the Bauer--Furuta invariant $\BF(X,\fraks)=\eta^2$ in $\pi_2^{s}$ holds.
Then the boundary Dehn twist $\delta_{\mathring{X}}$ is relatively exotic.    
\end{thm}

Here $\BF(X,\fraks)$ denotes the nonequivariant Bauer--Furuta invariant.
See Subsection~\ref{subsection: Families Bauer--Furuta invariants} for details.
This theorem yields, for example:

\begin{ex}[See Theorem~\ref{thm: elliptic boundary twist explicit}]
\label{ex: elliptic boundary twist explicit}  
Let $k\geq0$ and $l\geq1$, and let $p,q,r,s$ be positive odd integers with
$\gcd(p,q)=\gcd(r,s)=1$.
Assume that $(r,s)\neq(1,1)$ and that either $k$ is even or $\max\{p,q\}\geq5$.
Then the boundary Dehn twist on
\[
\left(E(4k+2)_{p,q}\#E(4l)_{r,s}\right)^\circ
\]
is relatively exotic.
\end{ex}

To the best of our knowledge, apart from the results of \cite[Theorem~1.2]{Lin20Dehn} for $(K3\#(S^2\times S^2))^\circ$ and of \cite{tilton2025boundarydehntwistpunctured} for $(K3\#K3)^\circ$, all previous results detecting nontrivial boundary Dehn twists on punctured 4-manifolds $\mathring{X}$ require $X$ to have a nonvanishing Seiberg--Witten invariant (\cite{BK19}, \cite[Proposition~1.2]{KM20Dehn}, \cite[Theorem~1.5]{baraglia2024irreducible4manifoldsadmitexotic}).
Consequently, these results do not apply to the examples in Example~\ref{ex: elliptic boundary twist explicit}, and our proof of Theorem~\ref{thm: eg puncture boundary Dehn twist} makes essential use of the Bauer--Furuta invariant.

\subsection{Seifert fibered Dehn twists on elliptic surfaces}

We now turn to Dehn twists along Seifert fibered 3-manifolds, which have been actively studied in recent years~\cite{konno2023exotic,KangParkTaniguchi,KLMME,OrsonPowell2022,miyazawaDehn,KLMME2,konno2024localizinggroupsexoticdiffeomorphisms,konno2025constraintslefschetzfibrationsfourdimensional,konno2026exoticdiffeomorphismsclosedaspherical}.

Following the first construction of exotic diffeomorphisms of irreducible 4-manifolds by Baraglia and the author~\cite{baraglia2024irreducible4manifoldsadmitexotic}, Lin, Mukherjee, Mu\~noz-Ech\'aniz, and the author recently constructed explicit examples of such diffeomorphisms arising as Dehn twists along Seifert fibered 3-manifolds~\cite{konno2025constraintslefschetzfibrationsfourdimensional}.
The main examples in \cite{konno2025constraintslefschetzfibrationsfourdimensional} are the Dehn twists along $\Sigma(2,3,7)$, which is the boundary of the Milnor fiber $M(2,3,7)$, embedded into $E(4k)_{p,q}$, where $p, q > 0$ are relatively prime odd integers outside a finite exceptional set.
We extend the result of~\cite{konno2025constraintslefschetzfibrationsfourdimensional} in two independent directions. First, we detect Dehn twists along the boundaries of other embedded Milnor fibers. Second, we extend the range of $(p,q)$ to include all relatively prime positive odd pairs except $(1,1)$.
Namely, we prove:

\begin{thm}
\label{thm: Dehn twists on elliptic surfaces}
Let $k\geq 1$, let $p,q$ be relatively prime positive odd integers, and let $1\leq i\leq 2k$.
Let $\delta$ denote the Dehn twist along $\del M(2,3,12i-5)$ in $E(4k)_{p,q}$, where $M(2,3,12i-5)$ is embedded in the standard way.
If $(p,q)\neq(1,1)$, then $\delta$ is an exotic diffeomorphism.
The same conclusion holds for the Dehn twist along $\del M(2,3,12i-1)$.
\end{thm}

This theorem follows from a constraint on the action on homology of a ``square root" of the Dehn twist.
See Theorem~\ref{thm: square root obstruction} for the precise statement.

\subsection{Dehn twists on nuclei and the detection of exotic structures}

Using Theorem~\ref{thm: Dehn twists on elliptic surfaces}, we provide a new characterization of the standard Gompf nucleus among exotic nuclei in terms of boundary Dehn twists.

Let $N(n)_{p,q}$ denote the Gompf nucleus obtained from $N(n)=N(n)_{1,1}$ by logarithmic transformations of multiplicities $p$ and $q$ in its interior.
The boundaries and their collar neighborhoods are naturally identified for all $p,q$.
When $n$ is even and $p,q$ are relatively prime positive odd integers, $N(n)_{p,q}$ is homeomorphic to $N(n)$ relative to this boundary identification \cite[Corollary (0.9), p. 334]{Boyer_with_bounadry}.
On the other hand, for $n\geq2$, $N(n)_{p,q}$ are pairwise nondiffeomorphic for distinct unordered pairs $\{p,q\}$ \cite[Corollary~8.3.14]{GS99}.
We detect the difference in smooth structures between $N(n)$ and $N(n)_{p,q}$ for $(p,q)\neq(1,1)$ using the boundary Dehn twist:

\begin{cor}
\label{cor: nuclei}
Let $k\geq 1$ be an integer, and let $p,q$ be relatively prime positive odd integers.
Then the boundary Dehn twist $\delta$ on $N=N(4k)_{p,q}$ is trivial in $\pi_0(\Diff_\del(N))$ if and only if $(p,q)=(1,1)$ (i.e. $N(4k)_{p,q}=N(4k)$).
\end{cor}

\begin{proof}
If $(p,q)\neq (1,1)$, Theorem~\ref{thm: Dehn twists on elliptic surfaces} implies that $\delta$ is nontrivial in $\pi_0(\Diff_\del(N))$.
Indeed, taking $i=2k$ in Theorem~\ref{thm: Dehn twists on elliptic surfaces}, we have that the Dehn twist along $\del M(2,3,24k-1)$ is nontrivial on $E(4k)_{p,q}$. Since $E(4k)_{p,q}$ splits into $E(4k)_{p,q} = M(2,3,24k-1) \cup N(4k)_{p,q}$, the Dehn twist along $\del M(2,3,24k-1)$ is the same as the twist along $\del N(4k)_{p,q}$, up to inversion.
Conversely, if $(p,q)=(1,1)$, then $\delta$ is trivial in $\pi_0(\Diff_\del(N))$ by \cite[Proposition 2.25]{KLMME}, which shows that the standard circle action on $\del N(4k)$ extends to a smooth circle action on $N(4k)$.
\end{proof}

%In this way, exotic diffeomorphisms ``detect" exotic smooth structures on 4-manifolds", 

%Corollary~\ref{cor: nuclei} means that the (non)triviality of the boundary Dehn twist is sensitive enough to the choice of smooth structure, in the sense that the triviality of the boundary Dehn twist characterizes the standard nucleus within exotic nuclei $N(n)_{p,q}$, provided that $n\equiv0 \bmod4$.

Corollary~\ref{cor: nuclei} shows that the (non)triviality of the boundary Dehn twist is sufficiently sensitive to the choice of smooth structure to distinguish the standard nucleus from exotic nuclei: for $n\equiv0\bmod4$, the triviality of the boundary Dehn twist characterizes the standard nucleus $N(n)$ among the nuclei $N(n)_{p,q}$ homeomorphic to it.

\subsection{Methods}

The results stated above are proved using the families Bauer--Furuta invariant.
To regard the families Bauer--Furuta invariant as an element of a stable homotopy group, we need to choose stable framings of the families Dirac index and the family of $H^+$ associated with a family of 4-manifolds.
For example, for a family of spin 4-manifolds over $S^1$, the quaternionic structure of the families Dirac index provides a canonical stable framing.
Kronheimer and Mrowka used this framing to detect the Dehn twist on $K3\#K3$~\cite{KM20Dehn}.
A major limitation of this approach has been the lack of a natural choice of framing without restrictive assumptions.
In particular, this is a key obstacle to extending the method of~\cite{KM20Dehn} to more general situations.

In this paper, we relax the requirement that the spin$^c$ structure arise from a spin structure and construct a framing of the family Dirac index for certain spin$^c$ structures $\fraks$ (Section~\ref{section A framing of the families Dirac index}).
The key assumption for constructing this framing is that $c_1(\fraks)$ is divisible by 8, as in the theorems stated above.
Assuming further that $X$ is a closed simply-connected spin 4-manifold with $\sigma(X)$ divisible by 32, we show that the families Bauer--Furuta invariant defined using this framing and the standard framing of $H^+$ gives rise to a group homomorphism
\[
\FBF_{\fraks}^{\divi} : \pi_0(\TDiff(X)) \to \pi_{d(\fraks)+2}^{s},
\]
where $\pi_0(\TDiff(X))$ denotes the Torelli group (Corollary~\ref{cor: FBF hom simply connected}).
Much of this paper is devoted to computing this version of the families Bauer--Furuta invariant for the Dehn twists discussed above, and the main nontriviality results stated above follow from these computations.

%Much of this paper is devoted to computing the families Bauer--Furuta invariant with these framings.
%All the results stated above are deduced from the nontriviality of these invariants.

\subsection*{Acknowledgements}

The author would like to thank two groups of collaborators: Abhishek Mallick and Masaki Taniguchi, and Jianfeng Lin, Anubhav Mukherjee, and Juan Mu\~noz-Ech\'aniz.
The possibility of a result such as Corollary~\ref{cor: nuclei} was discussed during these collaborations.
The author would also like to thank David Baraglia for his comments on an earlier version of this paper.
The author is grateful to the Okinawa Institute of Science and Technology (OIST) for its hospitality during the TSVP thematic program {\it New Frontiers in Gauge Theory, Topology, and Physics (TP26GT)}, during which this paper was written. The author was partially supported by JSPS KAKENHI Grant Numbers 25K00908 and 26K22244.

\subsection*{AI disclosure}

In the early stages of this project, the author considered a variant of the constraint ``$c_{1}(\Dsla)=w_{2}(H^{+})$'' established by Baraglia and the author \cite[Corollary 1.3]{BK19}, which he hopes to discuss elsewhere.
While investigating applications of that variant to Dehn twists, the author asked ChatGPT~6 whether the divisibility assumption $32\mid c_{1}(\fraks)$, used in \cite{baraglia2024irreducible4manifoldsadmitexotic} to apply \cite[Corollary 1.3]{BK19}, could be weakened.
ChatGPT suggested that divisibility by 8 might suffice and, more strongly, made an observation that yields a framing of the families Dirac index, a statement essentially equivalent to Lemma~\ref{lem: 8-divisibility framing}.
This suggestion led the author to explore a different approach to detecting Dehn twists, namely, computing the families Bauer--Furuta invariant for this framing, rather than applying a variant of \cite[Corollary 1.3]{BK19}.
From this point, ChatGPT was also used to assist the author in refining the arguments and proofreading the manuscript's English. In particular, the numerical computations in the proof of Lemma~\ref{lem: elliptic 8-divisible basic classes} were carried out with substantial assistance from ChatGPT.

\section{A framing of the families Dirac index}

\label{section A framing of the families Dirac index}

For a smooth family of 4-manifolds whose vertical tangent bundle is equipped with a spin structure, the families Dirac index naturally carries a quaternionic structure. When the base has dimension at most two, this structure canonically determines a stable framing of the families Dirac index.
This played an important role in the detection of Dehn twists by Kronheimer--Mrowka~\cite{KM20Dehn}, Lin~\cite{Lin20Dehn} and Miyazawa~\cite{miyazawaDehn}.
The main purpose of this section is to weaken the spin assumption and obtain a sufficient condition for obtaining such a canonical stable framing.
In short, when the 4-manifold is simply-connected, the divisibility of the first Chern class of a spin$^c$ structure $\fraks$ by 8 allows us to construct a stable framing of the families Dirac index, canonical up to homotopy, after choosing extensions of the underlying spin structure $\frakt$ and the line bundle $L$ to the family, where $\fraks=\frakt\otimes L^4$.

First, we note the following:

\begin{lem}
\label{lem:spin induces framing}
Let $B$ be a CW complex with $\dim B \leq 2$, and let $D \in KO(B)$.
Then a spin structure on $D$ induces a stable framing of the underlying real virtual bundle of $D$, unique up to homotopy.
\end{lem}

\begin{proof}
Write $D=V-\underline{\R}^N$, where $V \to B$ is an actual real vector bundle.
We may assume $\rank V \geq 3$ by stabilization.
Put $n=\rank V$.
The given spin structure on $D$, together with the standard spin structure on $\underline{\R}^N$, determines a spin structure on $V$.
Since $\Spin(n)$ is 1-connected if $n \geq 3$ and $\dim B \leq 2$, the principal $\Spin(n)$-bundle associated with this spin structure is trivial.
A trivialization of this principal bundle induces, via the canonical surjection $\Spin(n) \to \SO(n)$, a framing of $V$, and hence a stable framing of $D$.

Moreover, any two such trivializations differ by a map $B \to \Spin(n)$, which is null-homotopic since $\Spin(n)$ is 2-connected if $n\geq3$ and $\dim B \leq 2$.
Thus, the resulting stable framing is unique up to homotopy among those compatible with the given spin structure.
This completes the proof.
\end{proof}

In light of this lemma, we seek to construct a spin structure on the families Dirac index.
On the other hand, we have:

\begin{lem}[{Sasahira~\cite[Lemma 3.9]{MR2284407}}]
\label{lem:Sasahira}
Let $V \to B$ be a complex vector bundle.
Then a complex line bundle $R \to B$, together with an isomorphism $R^{\otimes 2} \cong \det_\C(V)$, canonically determines a spin structure on the underlying real bundle of $V$.
\end{lem}

Thus, our task is to construct an explicit ``square root" $R$ of the determinant line bundle of the families Dirac index $D$, together with an isomorphism identifying $R^{\otimes 2}$ with that determinant line bundle.

In what follows, throughout this section, let $X$ be an oriented closed spin smooth 4-manifold.
Let $\frakt$ be a spin structure on $X$, and let $L \to X$ be a complex line bundle.
Let $\fraks$ be the spin$^c$ structure on $X$ defined by $\fraks=\frakt \otimes L^4$, where $L^4$ denotes $L^{\otimes 4}$.
Note that $c_1(\fraks)=8c_1(L)$.
The reason for choosing the fourth power of $L$ shall be explained in Remark~\ref{rem: why8}.

\begin{rem}
\label{rem: simplyconn 8 divi}
In our applications, $X$ is simply-connected, so $\frakt$ is the unique spin structure on $X$. In this case, any spin$^c$ structure $\fraks$ whose first Chern class is divisible by $8$ can be written in the form $\fraks=\frakt \otimes L^4$ for some complex line bundle $L$ over $X$.
\end{rem}

Let $E \to B$ be an oriented smooth fiber bundle with fiber $X$.
Throughout this section, the base $B$ is a finite CW complex of $\dim B \leq 2$.
Assume that the spin structure $\frakt$ extends to a spin structure $\frakt_E$ on the vertical tangent bundle $T(E/B)$ of $E$, whose restriction to each fiber is isomorphic to $\frakt$.
Assume further that the complex line bundle $L$ extends to a complex line bundle $L_E \to E$ whose restriction to each fiber is isomorphic to $L$.
These extensions $\frakt_E$ and $L_E$ give a spin$^c$ structure $\fraks_E=\frakt_E \otimes L_E^4$ on $T(E/B)$ whose restriction to each fiber is isomorphic to $\fraks$.
We denote by $\ind_{\C}(\Dsla_{E,\fraks_E}) \in K(B)$ the index of the associated family of spin$^c$ Dirac operators for this family of spin$^c$ 4-manifolds with fiber $(X,\fraks)$.

We now give an explicit construction of a square root of $\det\nolimits_{\C}(\ind_{\C}(\Dsla_{E,\fraks_E}))$.
Fix a Hermitian metric and a unitary connection on $L$.
For each integer $j$, let
\[
\Dsla_{j} : \Gamma(X;S^+ \otimes L^j) \to \Gamma(X;S^-\otimes L^j)
\]
denote the spin Dirac operator for $\frakt$ twisted by $L^j$, i.e.,  the spin$^c$ Dirac operator for $\frakt \otimes L^j$.
We suppress the dependence on the choice of connection in our notation.
Note that, in particular, $\Dsla_{4}=\Dsla_{\fraks}$ is the spin$^c$ Dirac operator associated with $\fraks$.
Using $\frakt_E$, $L_E$, and families of Hermitian metrics and unitary connections, we obtain a family of twisted Dirac operators along the fibers of $E$ whose restriction to each fiber is the corresponding twisted Dirac operator.
We denote this family of twisted Dirac operators by
\[
\Dsla_{E,j} : \Gamma_B(E;S^+_E \otimes L_E^j) \to \Gamma_B(E;S^-_E \otimes L_E^j),
\]
where, for a vector bundle $W \to E$, $\Gamma_B(E;W)$ denotes the family of spaces of smooth sections of $W \to E$ along the fibers of $E \to B$.
Define complex line bundles $\lambda_{E,j} \to B$ by
\[
\lambda_{E,j} := \det\nolimits_{\C}(\ind_{\C}(\Dsla_{E,j})).
\]
Note that
\begin{align}
\label{eq: lambdaE4}
\lambda_{E,4} = \det\nolimits_{\C}(\ind_{\C}(\Dsla_{E,\fraks_E})).
\end{align}
Now, define a complex line bundle $R_E \to B$ by
\[
R_E = \lambda_{E,3}^2 \otimes \lambda_{E,2}^{-3} \otimes \lambda_{E,1}^2,
\]
where $\lambda_{E,j}^i=\lambda_{E,j}^{\otimes i}$.

We shall show that $R_E$ is indeed the desired square root of $\lambda_{E,4}$.
The existence of an isomorphism $R_E^{\otimes2}\cong\lambda_{E,4}$ can be seen
as follows.
Since $\dim E\leq6$, the bundle $L_E^{\oplus4}$ admits a nowhere-vanishing section, and hence splits off a trivial line bundle.
By the multiplicativity of the Euler class $\lambda_{-1}$ in complex $K$-theory, this implies $([L_E]-1)^4 =\lambda_{-1}(L_E^{\oplus4})=0$ in $K(E)$.
Expanding this relation, applying the families Dirac index, and taking determinant lines gives
\[
\lambda_{E,4}\cong
\lambda_{E,1}^4\otimes\lambda_{E,2}^{-6}
\otimes\lambda_{E,3}^4
=R_E^{\otimes2},
\]
where we use the canonical trivialization of $\lambda_{E,0}$ arising from the quaternionic structure.

However, to obtain a canonical framing, we need an isomorphism $R_E^{\otimes2}\cong\lambda_{E,4}$ that is canonical up to homotopy, rather than merely knowing that such an isomorphism exists.
The following lemma gives an explicit isomorphism:

\begin{lem}
\label{lem: 8-divisibility framing}
There exists an explicit isomorphism
\[
\phi : R_E^{\otimes 2} \to \det\nolimits_{\C}(\ind_{\C}(\Dsla_{E,\fraks_E}))
\]
whose homotopy class is uniquely determined by $(\frakt_E,L_E)$.
\end{lem}

\begin{proof}
First, note that the vector bundle $L_E^{\oplus4} \to E$ admits a nowhere-vanishing section, unique up to homotopy.
Indeed, the fibers of the sphere bundle $S(L_E^{\oplus4}) \to E$ are homeomorphic to $S^7$, which is 6-connected, and $\dim E \leq 6$.
Let $v : E \to L_E^{\oplus4}$ be a nowhere-vanishing section.

Fix a Hermitian metric on $L_E$, and equip $L_E^{\oplus 4}$ with the induced metric.
Let $v^\ast$ denote the metric dual of $v$.
Using exterior multiplication by $v$ and contraction with $v^\ast$, we obtain an isomorphism
\[
Q_v := v\wedge \bullet + v^\ast\lrcorner \bullet : \wedge^{\even} (L_E^{\oplus 4}) \to \wedge^{\odd} (L_E^{\oplus 4})
\]
of complex vector bundles, where the exterior powers are taken over $\C$.

As in the definition of $\Dsla_{E,j}$, given a Hermitian vector bundle $W \to E$ equipped with a unitary connection, let 
\[
\Dsla_{E,W} : \Gamma_B(E;S^+_E \otimes W) \to \Gamma_B(E;S^-_E \otimes W)
\]
denote by the family of spin Dirac operators twisted by $W$.
Consider the families of operators 
\begin{align*}
&\Dsla_E^{\even} := \Dsla_{E,\wedge^{\even} L_E^{\oplus4}}: \Gamma_B(E;S_E^+ \otimes \wedge^{\even} L_E^{\oplus4})
\to\Gamma_B(E;S_E^- \otimes \wedge^{\even} L_E^{\oplus4}),\\
&\Dsla_E^{\odd} := \Dsla_{E,\wedge^{\odd} L_E^{\oplus4}} : \Gamma_B(E;S_E^+ \otimes \wedge^{\odd} L_E^{\oplus4})
\to\Gamma_B(E;S_E^- \otimes \wedge^{\odd} L_E^{\oplus4}).
\end{align*}
%For brevity, write $\Gamma_B^\pm((L_E^j)^{\oplus i})$ for $\Gamma_B(E;S^\pm_E \otimes (L_E^j)^{\oplus i})$.
%We have operators
%\[
%\begin{aligned}
%\Dsla_E^{\even} := \Dsla_{E,0} \oplus \Dsla_{E,2}^{\oplus 6}\oplus \Dsla_{E,4}
%&:\Gamma_B^+(L^0_E)\oplus \Gamma_B^+((L_E^2)^{\oplus 6})\oplus \Gamma_B^+(L_E^4)\\
%&\to
%\Gamma_B^-(L^0_E)\oplus \Gamma_B^-((L_E^2)^{\oplus 6})\oplus \Gamma_B^-(L_E^4)
%\end{aligned}
%\]
%and
%\[
%\begin{aligned}
%\Dsla_E^{\odd} := \Dsla_{E,1}^{\oplus 4}\oplus \Dsla_{E,3}^{\oplus4}
%&:\Gamma_B^+(L_E^{\oplus 4})\oplus \Gamma_B^+((L_E^3)^{\oplus 4})\\
%&\to
%\Gamma_B^-(L_E^{\oplus 4}) \oplus \Gamma_B^+((L_E^3)^{\oplus 4})
%\end{aligned}
%\]
The isomorphism $Q_v$ induces isomorphisms between the domains and between the codomains of $\Dsla_E^{\even}$ and $\Dsla_E^{\odd}$, yielding the following diagram
\[
  \begin{CD}
     \Gamma_B(E;S_E^+ \otimes \wedge^{\even} L_E^{\oplus4}) @>{\Dsla_E^{\even}}>> \Gamma_B(E;S_E^- \otimes \wedge^{\even} L_E^{\oplus4}) \\
  @V{\id_{S_E^+}\otimes Q_v}VV    @V{\id_{S_E^-}\otimes Q_v}VV \\
     \Gamma_B(E;S_E^+ \otimes \wedge^{\odd} L_E^{\oplus4})   @>{\Dsla_E^{\odd}}>>  \Gamma_B(E;S_E^- \otimes \wedge^{\odd} L_E^{\oplus4}).
  \end{CD}
\]
This diagram commutes at the level of principal symbols.
More precisely, let $Q_v^\pm=\id_{S_E^\pm}\otimes Q_v$ and define
\[
\widetilde{\Dsla}_E^{\odd}
:=Q_v^-\circ\Dsla_E^{\even}\circ(Q_v^+)^{-1}.
\]
Then $\widetilde{\Dsla}_E^{\odd}$ and $\Dsla_E^{\odd}$ have the same family of principal symbols.
Indeed, the principal symbol of a twisted spin Dirac operator is given by Clifford multiplication on the spinor factor tensored with the identity on the twisting bundle, while $Q_v^\pm$ acts only on the twisting bundle.

Thus, $A=\Dsla_E^{\odd}-\tilde{\Dsla}_E^{\odd}$ is a family of zeroth-order operators. Consequently,
\begin{align}
\label{eq: Fred homotopy}
\tilde{\Dsla}_E^{\odd} + tA
\end{align}
for $t \in [0,1]$ defines a homotopy of families of Fredholm operators from $\tilde{\Dsla}_E^{\odd}$ to $\Dsla_E^{\odd}$.
This homotopy induces an isomorphism
\begin{align}
\label{eq: key lem psi}
\psi : \det\nolimits_{\C}(\ind_{\C}(\tilde{\Dsla}_E^{\odd})) \to \det\nolimits_{\C}(\ind_{\C}(\Dsla_E^{\odd}))
\end{align}
which is uniquely determined up to homotopy.
On the other hand, $Q_v$ defines an isomorphism
\[
\det\nolimits_{\C}(\ind_{\C}(\tilde{\Dsla}_E^{\odd})) \cong \det\nolimits_{\C}(\ind_{\C}(\Dsla_E^{\even})).    
\]
This composed with \eqref{eq: key lem psi} gives an isomorpshim
\begin{align}
\label{eq: composed psi}
\det\nolimits_{\C}(\ind_{\C}(\Dsla_E^{\even}))
\cong \det\nolimits_{\C}(\ind_{\C}(\Dsla_E^{\odd})).
\end{align}

Moreover, using the canonical isomorphisms
$\wedge^j(L_E^{\oplus n}) \cong (L_E^j)^{\oplus \tbinom{n}{j}}$,
we obtain
\begin{align*}
\wedge^{\even}(L_E^{\oplus4})
&\cong L_E^0 \oplus (L_E^2)^{\oplus6}\oplus L_E^4,\\
\wedge^{\odd}(L_E^{\oplus4})
&\cong L_E^{\oplus 4} \oplus (L_E^3)^{\oplus 4}.
\end{align*}
These decompositions induce canonical isomorphisms
\[
\begin{aligned}
\det\nolimits_{\C}(\ind_{\C}(\Dsla_E^{\even}))
&\cong \lambda_{E,0} \otimes \lambda_{E,2}^{\otimes 6} \otimes \lambda_{E,4} \cong \lambda_{E,2}^{\otimes 6} \otimes \lambda_{E,4},\\
\det\nolimits_{\C}(\ind_{\C}(\Dsla_E^{\odd}))
&\cong \lambda_{E,1}^{\otimes 4} \otimes \lambda_{E,3}^{\otimes 4}.
\end{aligned}
\]
Here we use the canonical trivialization of $\lambda_{E,0}$ arising from the quaternionic structure on the family of spin Dirac operators.
These isomorphisms, together with the isomorphism \eqref{eq: composed psi}, yield an isomorphism
\[
\phi : (\lambda_{E,3}^{\otimes 2} \otimes \lambda_{E,2}^{\otimes(-3)} \otimes \lambda_{E,1}^{\otimes 2})^{\otimes 2} \to \lambda_{E,4}.
\]

Recall that any two nowhere-vanishing sections $v$ are homotopic. Moreover, the spaces of other auxiliary choices, namely, fiberwise Riemannian metrics on $E$, Hermitian metrics and unitary connections on $L_E$, are contractible.
Thus, the homotopy class of $\phi$ is independent of these auxiliary choices and is canonically determined by $(\frakt_E,L_E)$.
Recalling \eqref{eq: lambdaE4}, we obtain the desired isomorphism
$\phi : R_E^{\otimes 2} \to \det\nolimits_{\C}(\ind_{\C}(\Dsla_{E,\fraks_E}))$.
\end{proof}

By construction, the pair $(R_E,\phi)$ is natural under pullback and under isomorphisms of families preserving the pair $(\frakt_E,L_E)$.

\begin{rem}
\label{rem: why8}
One reason for choosing the fourth power of $L$, corresponding to the divisibility of $c_1(\fraks)$ by $8$, is the following.
In the above argument, the uniqueness up to homotopy of a nowhere-vanishing section $v$ is guaranteed by the fact that the fibers of the sphere bundle $S(L_E^{\oplus4})$ are homeomorphic to $S^7$, which is 6-connected. Even in our applications, where $B=S^1$ and hence $\dim E=5$, the sphere fiber must be at least $5$-connected to ensure uniqueness up to homotopy by this argument. This condition is not satisfied by $S(L_E^{\oplus m})$ for $m\leq3$.
Thus, the fourth power is the smallest power for which the above construction works under our assumptions.
\end{rem}

In summary, we obtain:

\begin{cor}
\label{cor: Dirac framing by 8}
A pair of extensions $(\frakt_E,L_E)$ of $(\frakt,L)$ to $E$ canonically determines, up to homotopy, a stable framing of the underlying real virtual bundle of $\ind_{\C}(\Dsla_{E,\fraks_E})$, where $\fraks_E=\frakt_E \otimes L_E^4$.
\end{cor}

\begin{proof}
This follows from Lemma~\ref{lem:spin induces framing}, Lemma~\ref{lem:Sasahira}, and Lemma~\ref{lem: 8-divisibility framing}.
\end{proof}

\begin{defi}
\label{defi:8framing}
We denote by $\xi_{\Dsla}^{\divi}$ the stable framing of the underlying real virtual bundle of $\ind_{\C}(\Dsla_{E,\fraks_E})$ constructed in Corollary~\ref{cor: Dirac framing by 8}, and call it the \textit{8-divisibility framing}.
When we wish to emphasize its dependence on $(\frakt_E,L_E)$, we write $\xi_{\Dsla}^{\divi}(\frakt_E,L_E)$.
By abuse of terminology, we also refer to the isomorphism $\phi : R_E^{\otimes 2} \to \det\nolimits_{\C}(\ind_{\C}(\Dsla_{E,\fraks_E}))$ given in Lemma~\ref{lem: 8-divisibility framing} as the 8-divisibility framing.
We call $(R_E,\phi)$ {\it square root data} for $E$.
\end{defi}

%\begin{rem}
%Put $x=L_E-1\in K(E)$. Since $\dim E \leq 6$, we have $x^4=0$, and hence
%\[
%L_E^4=1+4x+6x^2+4x^3.
%\]
%By the additivity of the families index, it follows that
%\[
%\ind_{\C}(\Dsla_{E,4}) = -\ind_{\C}(\Dsla_{E,0}) + 2\left(2\ind_{\C}(\Dsla_{E,1})-3\ind_{\C}(\Dsla_{E,2})+2\ind_{\C}(\Dsla_{E,3})\right).
%\]
%Since $\Dsla_{E,0}$ is a family of spin Dirac operators, its index carries a quaternionic structure and hence has vanishing first Chern class.Since the remaining terms are divisible by $2$, we conclude that $c_1(\ind_{\C}(\Dsla_{E,4}))$ is divisible by $2$.
%This implies that $\det_{\C}$
%\end{rem}

\section{Computing families Bauer--Furuta invariants}

\subsection{Families Bauer--Furuta invariants}
\label{subsection: Families Bauer--Furuta invariants}

We now briefly review the families Bauer--Furuta invariants. For further details, see, for example, \cite{BF04,BK19,KM20Dehn,konno2025constraintslefschetzfibrationsfourdimensional}.
Let $(X,\fraks)$ be a closed smooth spin$^c$ 4-manifold with $b_1(X)=0$.
Let $E \to B$ be a smooth family of spin$^c$ 4-manifolds with fiber $(X,\fraks)$.
That is, $E \to B$ is a smooth fiber bundle with fiber $X$ whose vertical tangent bundle $T(E/B)$ is equipped with a spin$^c$ structure $\fraks_E$ whose restriction to each fiber is isomorphic to $\fraks$.
After fixing a fiberwise Riemannian metric on $E$, we can consider the family of Seiberg--Witten equations along the fibers of $E$.
By taking a finite-dimensional approximation of the family of Seiberg--Witten equations, we obtain a pointed continuous map between Thom spaces
\begin{align}
\label{eq: fin dim approx Psi}
\Psi_E : \Th(V^+ \oplus U^+) \to \Th(V^- \oplus U^-),
\end{align}
where $V^+, V^- \to B$ are finite-dimensional complex vector bundles satisfying $V^+ - V^- = \ind_{\C}(\Dsla_{E,\fraks_E})$ in $K(B)$, and $U^+, U^- \to B$ are finite-dimensional real vector bundles satisfying $U^- - U^+ = H^+(E)$ in $KO(B)$.
Here $H^+(E) \to B$ is a real vector bundle of rank $b^+(X)$ whose isomorphism class is canonically determined by $E$. For a choice of fiberwise Riemannian metric on $E$, it can be realized as the bundle of self-dual harmonic 2-forms along the fibers of $E$.

Let $S^1$ act on $V^\pm$ by scalar multiplication and trivially on $U^\pm$. The map $\Psi_E$ is equivariant with respect to the induced $S^1$-actions on the Thom spaces.

When $B=\{\pt\}$, we write $\Psi_X$ for $\Psi_E$, and the construction reduces to the usual (i.e., unparameterized) nonequivariant Bauer--Furuta construction for $(X,\fraks)$.
To obtain this invariant as an element of $\pi_{d(\fraks)+1}^{s}$, we fix an orientation of $H^+(X)$. This determines a homotopy class of orientation-preserving isomorphisms from $H^+(X)$ to $\R^{b^+(X)}$. Similarly, the complex structure on $\ind_{\C}\Dsla_{\fraks}$ induces a canonical orientation on its realification $\ind_{\R}\Dsla_{\fraks}$, which has virtual dimension $(c_1(\fraks)^2-\sigma(X))/4$. This determines a stable framing of $\ind_{\R}\Dsla_{\fraks}$ uniquely up to homotopy.
By forgetting the $S^1$-equivariance of $\Psi_X$ and using these framings, we obtain a nonequivariant stable map
\[
\BF(X,\fraks) : S^{d(\fraks)+1} \to S^0.
\]
We call its stable homotopy class the {\it nonequivariant Bauer--Furuta invariant} of $(X,\fraks)$.
Here $d(\fraks)$ denotes the Seiberg--Witten formal dimension, given by
\[
d(\fraks)=\frac{c_1(\fraks)^2-2\chi(X)-3\sigma(X)}{4}.
\]
We obtain an element of the stable homotopy group $\pi_{d(\fraks)+1}^{s}$ from this stable map, which we also denote by $\BF(X,\fraks) \in \pi_{d(\fraks)+1}^{s}$.

We shall use the following fact about the unparameterized Bauer--Furuta invariant:

\begin{lem}[{Bauer~\cite[Proposition~4.4]{B04}}]
\label{lem: BF equals eta iff SW odd}
Let $X$ be a closed, oriented, smooth spin 4-manifold with $b_1(X)=0$, and let $\fraks$ be a spin$^c$ structure on $X$ such that $c_1(\fraks)$ is divisible by $8$.
Suppose that $d(\fraks)=0$.
Then $b^+(X)\equiv 3 \bmod 4$, and $\BF(X,\fraks)=\eta$ if and only if the Seiberg--Witten invariant $\SW(X,\fraks)$ is odd.
\end{lem}

\begin{proof}
As seen in the previous section, the divisibility of $c_1(\fraks)$ by 8 implies that $r=\ind_{\C}\Dsla_{\fraks}$ is even.
It follows from this combined with $d(\fraks)=0$ that $b^+(X)=2r-1\equiv3$ mod 4.
Moreover, $d(\fraks)=0$ implies that $\fraks$ is induced by an almost complex structure on $X$
(cf.~\cite[p.36, the paragraph following Corollary~4.2]{B04}).
The assertion now follows from \cite[Proposition~4.4]{B04}.
\end{proof}

In the parameterized setup, our main focus is the case $B=S^1$.
For every family $E \to S^1$ considered in this paper, we assume that $H^+(E)$ is orientable and fix an orientation compatible with the orientation of $H^+(X)$ chosen above.
Assume that the underlying real virtual bundle of $\ind_{\C}(\Dsla_{E,\fraks_E})$ and the real vector bundle $H^+(E)$ are equipped with stable framings $\xi_{\Dsla}$ and $\xi_{H^+}$, respectively.
Throughout, we require these framings to be compatible with the canonical orientation induced by the complex structure on the Dirac index and the chosen orientation of $H^+(E)$, respectively.
Using these framings, we obtain a stable map
\begin{align}
\label{eq:tilde FBF}
\widetilde{\FBF}(E,\fraks_E,\xi_{\Dsla},\xi_{H^+}) : S^{d(\fraks)+1} \wedge B_+ \to S^0
\end{align}
in the nonequivariant stable homotopy category.
Choose a point $x\in B=S^1$. The natural inclusion $i_x:S^0\hookrightarrow B_+$ and the natural projection $p:B_+\to S^0$ satisfy $p\circ i_x=\id_{S^0}$. Since the cofiber of $i_x$ is $(B,x)$, the associated cofiber sequence splits stably, giving $B_+\simeq S^0\vee S^1$.
Consequently, the domain of \eqref{eq:tilde FBF} admits a stable splitting
\begin{align}
\label{eq:stable spliting}
S^{d(\fraks)+1}\wedge B_+
\simeq 
S^{d(\fraks)+1}\vee S^{d(\fraks)+2}.
\end{align}
We define
\[
\FBF(E,\fraks_E,\xi_{\Dsla},\xi_{H^+})\in\pi_{d(\fraks)+2}^{s}
\]
to be the component of \eqref{eq:tilde FBF} corresponding to the $S^{d(\fraks)+2}$-summand.
We call this element the {\it nonequivariant families Bauer--Furuta invariant} of $(E,\fraks_E,\xi_{\Dsla},\xi_{H^+})$.
Note that the component of \eqref{eq:tilde FBF} corresponding to the $S^{d(\fraks)+1}$-summand is the unparameterized Bauer--Furuta invariant $\BF(X,\fraks)\in\pi_{d(\fraks)+1}^{s}$.

Let $f : X \to X$ be an orientation-preserving diffeomorphism satisfying $f^\ast \fraks \cong \fraks$.
Let $E_f \to S^1$ denote the mapping torus of $f$.
A choice of lift $\tilde{f}$ of $f$ to an automorphism of the spin$^c$ structure $\fraks$ induces a spin$^c$ structure $\fraks_{E_f}$ on the vertical tangent bundle of $E_f$, whose restriction to each fiber is isomorphic to $\fraks$.
Any two such lifts differ by a gauge transformation $X \to S^1$.
Since $b_1(X)=0$, the gauge group $\Map(X,S^1)$ is path-connected, and hence the lift is unique up to homotopy.
We therefore suppress the choice of $\tilde{f}$ from our notation and define
\[
\FBF(X,\fraks,f,\xi_{\Dsla},\xi_{H^+})
:=\FBF(E_f,\fraks_{E_f},\xi_{\Dsla},\xi_{H^+}) \in \pi_{d(\fraks)+2}^{s}.
\]

Moreover, if $f$ acts trivially on $H^2(X;\R)$, the vector bundle over $S^1$ associated with $E_f$, with fiber $H^2(X;\R)$ has a canonical trivialization.
Since the space of maximal positive-definite subspaces of $H^2(X;\R)$ is contractible, a choice of reference framing determines a framing of $H^+(E_f)$, unique up to homotopy.
We denote this framing by $\xi_{H^+}^{\std}$ and call it the {\it standard framing} of $H^+(E_f)$.

Since all stable framings are assumed to be compatible with the fixed orientations, the difference between two stable framings $\xi_{\Dsla}$ and $\xi_{\Dsla}'$ of the real families Dirac index over $S^1$ determines an element of $\pi_1(\SO)\cong\Z/2$. The same holds for two stable framings of $H^+(E_f)$.
The following formula describes how the families Bauer--Furuta invariant changes under a change of framing. It is a slight generalization of \cite[Proposition~3.6]{konno2025constraintslefschetzfibrationsfourdimensional}, allowing $d(\fraks)$ to be nonzero. The formula is also essentially equivalent to \cite[Proposition~2.1]{baraglia2026exoticdiffeomorphismsreducible4manifolds} due to Baraglia. However, passing from the formulation in \cite{baraglia2026exoticdiffeomorphismsreducible4manifolds} to ours requires some translation, so we include a proof in our notation.
Let $\eta \in \pi_1^{s}\cong\Z/2$ denote the generator.

\begin{lem}[{Baraglia \cite[Proposition 2.1]{baraglia2026exoticdiffeomorphismsreducible4manifolds}}]
\label{lem:framing_change}
Let $f : X \to X$ be an orientation-preserving diffeomorphism satisfying $f^\ast\fraks\cong\fraks$, and let $E_f\to S^1$ be its mapping torus equipped with a family spin$^c$ structure extending $\fraks$.
Fix stable framings $\xi_{\Dsla}$ and $\xi_{H^+}$ of the real families Dirac index and $H^+(E_f)$, respectively.
Then, for any $a, b\in\Z/2$, we have
\[
\FBF(X,\fraks,f,\xi_{\Dsla}+a,\xi_{H^+}+b)
= \FBF(X,\fraks,f,\xi_{\Dsla},\xi_{H^+})
+ (a+b) \cdot \eta \cdot \BF(X,\fraks).
\]
\end{lem}

\begin{proof}
Put $k=d(\fraks)+1$.
Represent the stable map \eqref{eq:tilde FBF} by a map
\[
F : S^{N+k} \wedge S^1_+ \to S^N
\]
for large enough $N>0$.
Let $l:S^1\to\SO(N+k)$ be a based loop representing $a+b\in\pi_1(\SO)\cong\Z/2$. Changing the Dirac framing by $a$ and $b$ corresponds, after stabilization, to precomposing $F$ with the fiberwise self-map
\[
T:S^{N+k}\wedge S^1_+\to S^{N+k}\wedge S^1_+
\]
define by $T(v,t)=(l(t)v,t)$.

Recall that we used the stable splitting \eqref{eq:stable spliting}, represented by
\[
S^{N+k}\wedge S^1_+\simeq
S^{N+k}\vee S^{N+k+1},
\] 
to define the invariant $\FBF(E,\fraks_E,\xi_{\Dsla},\xi_{H^+})$.
Under this splitting, 
\[
T : S^{N+k}\vee S^{N+k+1} \to S^{N+k}\vee S^{N+k+1}
\]
is obviously of the form
\[
\begin{pmatrix}
1&A\\
0&1
\end{pmatrix}.
\]
By the definition of $T$, the component $A : S^{N+k+1} \to S^{N+k}$ is given by $A=J(a+b)$, where $J : \pi_1(\SO) \to \pi_{1}^{s}$ is the $J$-homomorphism.

On the other hand, the two components of $F$ under the same splitting \eqref{eq:stable spliting} are $\BF(X,\fraks)$ and $\FBF(X,\fraks,f,\xi_{\Dsla},\xi_{H^+})$, respectively. Therefore, precomposing with $T$ changes the second component by $J(a+b)\cdot\BF(X,\fraks)$. Since $J(a+b)=(a+b)\eta$, the desired formula follows.
\end{proof}

We also note the following naturality under the pull-back of the family along a map $S^1 \to S^1$.
Let $p : S^1 \to S^1$ be a continuous map with degree $k \in \Z$.
Then the pull-back of a mapping torus $E_f \to S^1$ under $p$ gives also a family $p^\ast E_f \to S^1$, which is isomorphic to the mapping torus $E_{f^k}$ of $f^k$.
Given stable framings $\xi_{\Dsla},\xi_{H^+}$ for $E_f$, we obtain the pull-back framings $p^\ast\xi_{\Dsla}, p^\ast\xi_{H^+}$ for the family $p^\ast E_f \to S^1$.

\begin{lem}
\label{lem: pullback}
In the setup of Lemma~\ref{lem:framing_change}, let $p : S^1 \to S^1$ be a continuous map with degree $k \in \Z$.
Then we have 
\[
\FBF(X,\fraks,f^k,p^\ast\xi_{\Dsla},p^\ast\xi_{H^+})
= k\cdot \FBF(X,\fraks,f,\xi_{\Dsla},\xi_{H^+})
\]
in $\pi_{d(\fraks)+2}^{s}$.
\end{lem}

\begin{proof}
This is obvious from the definition of the families Bauer--Furuta invariant.    
\end{proof}

\subsection{Twisted families}

Following \cite[Section 4]{KM20Dehn}, we compute the families Bauer--Furuta invariant of the family twisted by the deck transformation of the spin structure.
However, we allow spin$^c$ structures that are not spin and use the 8-divisibility framing (Definition~\ref{defi:8framing}) instead of the quaternionic framing employed in \cite[Section 4]{KM20Dehn}.

As in Section~\ref{section A framing of the families Dirac index}, let $X$ be an oriented closed spin smooth 4-manifold with $b_1(X)=0$, $\frakt$ be a spin structure on $X$, and let $L \to X$ be a complex line bundle.
Let $\fraks$ be the spin$^c$ structure on $X$ defined by $\fraks=\frakt \otimes L^4$.

For $j \geq 0$, recall that $\Dsla_j$ denotes the spin$^c$ Dirac operator for $\frakt\otimes L^j$.
Define $r_j =\ind_{\C}\Dsla_{j}$, and set $r:=r_4=\ind_{\C}\Dsla_{\fraks}$.
By the index theorem, we have
\begin{align}
\label{eq:rj compute}
r_j = \frac{j^2 c_1(L)^2}{2} - \frac{\sigma(X)}{8}.
\end{align}
In particular, $r$ is even, since $\sigma(X)/8$ is even by Rokhlin's theorem.

Let $\tau: \frakt \to \frakt$ denote the nontrivial deck transformation of the spin structure $\frakt$. It is the unique non-trivial element in the kernel of the natural homomorphism from the group of spin automorphisms of $(X,\frakt)$ to the diffeomorphism group of $X$.
The mapping torus $\frakt^\tau$ of $\tau$ is the families spin structure with fiber $\frakt$ on the trivialized bundle $X \times S^1 \to S^1$.
We consider the trivial extension $\underline{L}=L \times S^1$ of $L$. 

Again by Rokhlin's theorem, $\sigma(X)/16$ is an integer.
Define $\overline{\sigma}(X) \in \Z/2$ by 
\[
\overline{\sigma}(X) := \frac{\sigma(X)}{16} \bmod 2.
\]

\begin{prop}
\label{prop:twisted family}
We have
\[
\FBF(X,\fraks,\id_X,\xi_{\Dsla}^{\divi}(\frakt^\tau, \underline{L}),\xi_{H^+}^{\std})
= \overline{\sigma}(X) \cdot \eta \cdot \BF(X,\fraks)
\]
in $\pi^{s}_{d(\fraks)+2}$.
\end{prop}

\begin{proof}
The basic strategy is similar to that used in \cite[Proof of Proposition 4.1]{KM20Dehn}: we compare two framings of the families Dirac index.

Let
\[
\Psi_X : (V^+ \oplus U^+)^+ \to (V^- \oplus U^-)^+
\]
be a finite-dimensional approximation of the Seiberg--Witten map for $(X,\fraks)$.
Here $(-)^+$ denotes the one-point compactification.
A finite-dimensional approximation $\Psi_E$ of the families Seiberg--Witten map for the family $E=(X \times S^1,\frakt^{\tau}\otimes\underline{L}^4) \to S^1$ is obtained as the mapping torus of $\Psi_X$ with respect to the action induced by $\tau$.
More precisely, $\Psi_E$ is obtained from
\[
\Psi_X \times \id_{[0,1]} : (V^+ \oplus U^+)^+ \times [0,1] \to (V^- \oplus U^-)^+ \times [0,1]
\]
by identifying the fibers over $t=0$ and $t=1$ via multiplication by $-1$ on the spinor summands $V^\pm$ and the identity on $U^\pm$.

We now examine the action induced by this multiplication by $\tau$ on the fibers of $R_E \to S^1$ and $\lambda_{E,4} \to S^1$.
We write $R$ and $\lambda_j$ for the corresponding line bundles over a point; that is,
$R=R_{(X,\fraks)\times\{\pt\}}$ and $\lambda_j=\lambda_{j,(X,\fraks)\times\{\pt\}}$.
Then 
\[
\lambda_j = \wedge^{\max}\Ker(\Dsla_j) \otimes \wedge^{\max}(\Coker(\Dsla_j))^\ast.
\]
Since $\tau$ acts by multiplication by $-1$ on both $\Ker(\Dsla_j)$ and $\Coker(\Dsla_j)$, its induced action on $\lambda_j$ is multiplication by
\[
(-1)^{\dim_{\C}\Ker(\Dsla_j)}\cdot((-1)^{\dim_{\C}\Coker(\Dsla_j)})^{-1}
=(-1)^{r_j}.
\]
Since $r=r_4$ is even, as noted above, $\tau$ acts trivially on $\lambda_4$.
On the other hand, the action of $\tau$ on $R$ is given by multiplication by $(-1)^{2r_3-3r_2+2r_1}$.
By \eqref{eq:rj compute}, a straightforward computation gives
\[
2r_3-3r_2+2r_1
= 4c_1(L)^2-\frac{\sigma(X)}{8},
\]
which is even by Rokhlin's theorem.
Hence, $\tau$ acts trivially on $R$ as well.

Since $\tau$ acts trivially on both $R$ and $\lambda_4$, the complex line bundles $R_E \to S^1$ and $\lambda_{E,4} \to S^1$ are obtained by gluing the endpoints of $R \times [0,1]$ and $\lambda_4 \times [0,1]$, respectively, via the identity maps. Thus, both bundles admit natural product trivializations after fixing identifications $R\cong\C$ and $\lambda_4\cong\C$.

Moreover, we claim that the isomorphism $\phi : R_E^{\otimes2} \to \lambda_{E,4}$ constructed in Lemma~\ref{lem: 8-divisibility framing} is constant with respect to these product structures.

To see this, recall that, in the proof of Lemma~\ref{lem: 8-divisibility framing}, we chose a nowhere-vanishing section $v$ of $L_E^{\oplus4}$, together with Hermitian metrics and unitary connections on the relevant bundles.

Consider the product family $E'=(X,\frakt \otimes L^4)\times[0,1]$ over $[0,1]$, and choose all these auxiliary data to be independent of the parameter $t\in[0,1]$.
Then the associated operators $Q_v$, $\Dsla_{E'}^{\even}$, and $\Dsla_{E'}^{\odd}$ are also independent of $t$.

Recall that these operators, together with the homotopy \eqref{eq: Fred homotopy} between them used in the construction, determine an isomorphism \eqref{eq: key lem psi}, 
\[
\psi : \det\nolimits_{\C}(\ind_{\C}(\tilde{\Dsla}_{E'}^{\odd})) \to \det\nolimits_{\C}(\ind_{\C}(\Dsla_{E'}^{\odd})),
\]
of complex line bundles over $[0,1]$.
Since the operators and the homotopy between them are independent of $t \in [0,1]$, the isomorphism $\psi$ can also be chosen to be independent of $t$.

Since $\tau$ acts trivially on $L$ and by multiplication by $-1$ on the spinor bundles, it commutes with $Q_v$, $\Dsla_{E'}^{\even}$, and $\Dsla_{E'}^{\odd}$ at $t=1$.
It also commutes with the homotopy connecting the relevant operators, and thus $\psi$ commutes with the induced actions of $\tau$ on the domain and codomain of $\psi$.

Furthermore, $\tau$ acts trivially on $\lambda_0$ by \eqref{eq:rj compute} and Rokhlin's theorem and therefore preserves its canonical quaternionic trivialization.

Combining these observations, we conclude that the isomorphism $R_{E'}^{\otimes2}\to\lambda_{E',4}$ constructed in Lemma~\ref{lem: 8-divisibility framing} for $E'$ can be chosen to be independent of $t$ and is compatible with the action of $\tau$ at $t=1$.
Consequently, it descends to the isomorphism $\phi:R_E^{\otimes2}\to\lambda_{E,4}$ for the family $E=E'/\tau$, which is constant with respect to the product structures
$R_E^{\otimes2}= R^{\otimes2}\times S^1$, $\lambda_{E,4}=\lambda_4\times S^1$.
This completes the proof of the claim.

We have thus shown that the 8-divisibility framing of $\lambda_{E,4}$ agrees, up to homotopy, with its product trivialization. In other words, the 8-divisibility framing is given, up to homotopy, by
\[
\alpha\times\id_{S^1}:\lambda_{E,4}=\lambda_4\times S^1 \to \C\times S^1
\]
for some isomorphism $\alpha:\lambda_4 \to \C$.

We now compare this 8-divisibility framing with the framing induced by the half-rotation arising from the $S^1$-action on the Seiberg--Witten map, as used in \cite[Proof of Proposition 4.1]{KM20Dehn}.
Recall that the scalar $S^1$-action on the spinor summands $V^\pm$ of the finite-dimensional approximation provides a path of complex-linear automorphisms from $\id_{V^\pm}$ to the action of $\tau$, given by multiplication by $-1$.
This path induces trivializations of the bundles $(V^\pm\times[0,1])/\tau\to S^1$.
We refer to this framing as the half-rotation framing and denote it by $\xi_{\Dsla}^{\pi}$.
Since $\Psi_X$ is $S^1$-equivariant, the framing $\xi_{\Dsla}^{\pi}$ identifies $\Psi_E$ with the product map $\Psi_X\times\id_{S^1}$.
Thus, we have
\begin{align}
\label{eq: vanishing FBF}
\FBF(X,\fraks,\id_X,\xi_{\Dsla}^{\pi},\xi_{H^+}^{\std})=0
\end{align}
in $\pi_{d(\fraks)+2}^{s}$.

As in the above computation of the action of $\tau$ on $\lambda_j$, the half-rotation $[0,1] \to \GL_{\C}(V^\pm)$ given by $t \mapsto e^{i\pi t}$ induces a path $[0,1] \to \GL_{\C}(\lambda_4)$ given by $t \mapsto e^{i\pi rt}$.
Since $r$ is even, this path is a loop, with winding number $r/2$.
Thus, the half-rotation framing of $\lambda_{E,4}$ differs from the product trivialization $\lambda_{E,4}=\lambda_4\times S^1$, which agrees up to homotopy with the 8-divisibility framing as shown above, by a loop of winding number $r/2$.
After stabilization, the difference between the corresponding stable framings is
\[
\frac{r}{2} \bmod 2 \in \pi_1(\SO) \cong \Z/2.
\]
By \eqref{eq:rj compute}, this equals $\overline{\sigma}(X)=\sigma(X)/16\bmod 2$.
The assertion of the proposition follows from this combined with Lemma~\ref{lem:framing_change} and \eqref{eq: vanishing FBF}.
\end{proof}

\subsection{\texorpdfstring{$\FBF$}{FBF}
defined over the Torelli group}
\label{subsection FBF homomorphism}

Let $\TDiff(X)$ denote the Torelli diffeomorphism group,
\[
\TDiff(X) = \{f \in \Diff(X) \mid f^\ast=\id \text{ on }H^\ast(X;\Z) \}.
\]
The group of path-components $\pi_0(\TDiff(X))$ is called the Torelli group of $X$.
In this subsection, under the assumption that $X$ is spin, simply-connected, $\sigma(X)$ is divisible by 32, and $\fraks$ is a spin$^c$ structure on $X$ such that $c_1(\fraks)$ is divisible by 8, we construct a well-defined group homomorphism
\[
\FBF_{\fraks}^{\divi} : \pi_0(\TDiff(X)) \to \pi_{d(\fraks)+2}^{s}
\]   
using the families Bauer--Furuta invariant for the 8-divisibility framing and the standard framing for $H^+$.
In particular, this invariant can be used to detect nontrivial smooth mapping classes.

As a first step, we establish the following independence result, which is a variant of the proof of the above Proposition~\ref{prop:twisted family}:

\begin{prop}
\label{prop: 8 framing can detect nontriviality}
Let $X$ be a closed smooth spin 4-manifold with $b_1(X)=0$ and $\sigma(X)$ divisible by $32$.
Let $\frakt$ be a spin structure on $X$, and let $L\to X$ be a complex line bundle.
Let $f:X\to X$ be an orientation-preserving diffeomorphism satisfying $f^\ast\frakt\cong\frakt$, $f^\ast L\cong L$, and $f^\ast=\id$ on $H^2(X;\R)$.
Set $\fraks=\frakt\otimes L^4$.
Then the invariant
\begin{align}
\label{eq: FBF in triviality lemma}
\FBF(X,\fraks,f,\xi_{\Dsla}^{\divi}(\frakt_E,L_E),\xi_{H^+}^{\std})
\in \pi_{d(\fraks)+2}^{s}
\end{align}
is independent of the choices of extensions $\frakt_E$ and $L_E$ of $\frakt$ and $L$ to the mapping torus $E_f$.
\end{prop}

\begin{proof}
Put $E=E_f$.
First, the assumption $b_1(X)=0$ implies that an extension $L_E$ of $L$ to $E$ is unique up to isomorphism.
Indeed, any two choices of gluing isomorphisms $f^\ast L\cong L$ differ by an automorphism of $L$, represented by a map $X\to S^1$.
Since $H^1(X;\Z)=0$, this map is null-homotopic, and hence the resulting line bundles over $E$ are isomorphic.

Thus, we may fix an extension $L_E$ of $L$ and need only compare the two extensions of $\frakt$.
Fix a lift $\tilde{f}:\frakt\to\frakt$ of $f:X\to X$.
The other lift is $\tau\circ\tilde{f}$, where $\tau:\frakt\to\frakt$ is the unique nontrivial deck transformation.
These two lifts determine extensions $\frakt_E^+$ and $\frakt_E^-$ of $\frakt$ to $E$, respectively.
Set $\xi_{\Dsla}^\pm=\xi_{\Dsla}^{\divi}(\frakt_E^\pm,L_E)$.

Let $\Psi_E^+$ be an $S^1$-equivariant finite-dimensional approximation of the families Seiberg--Witten map associated with $(E,\frakt_E^+\otimes L_E^4)$.
Pulling back $\Psi_E^+$ along the quotient map $[0,1]\to S^1$, we obtain a finite-dimensional approximation $\Psi_{X\times [0,1]}$ for the pullback family over $[0,1]$, whose underlying manifold bundle is $X\times[0,1]$.
We recover $\Psi_E^+$ from $\Psi_{X\times [0,1]}$ by identifying the fibers over $t=0$ and $t=1$ via the gluing maps induced by $\tilde{f}$ and the chosen gluing isomorphism for $L_E$.

On the other hand, a finite-dimensional approximation $\Psi_E^-$ associated with $(E_f,\frakt_E^-\otimes L_E^4)$ can be obtained from the same $\Psi_{X\times [0,1]}$ by composing these gluing maps with the action of $\tau$.
Here $\tau$ acts by multiplication by $-1$ on the spinor summands and trivially on the remaining summands.

For each sign $\pm$, let $R_E^\pm\to S^1$ and $\lambda_{E,4}^\pm\to S^1$ denote the complex line bundles corresponding to $R_E, \lambda_{E,4}$ associated with $(E,\frakt_E^\pm,L_E)$.
Similarly, let $\phi^\pm:(R_E^\pm)^{\otimes2}\to\lambda_{E,4}^\pm$ denote the corresponding 8-divisibility framings.
As in the proof of Proposition~\ref{prop:twisted family}, $\tau$ acts trivially on $R$ and $\lambda_4$. Thus, by the same argument used to establish the existence of product trivializations of $R_E$ and $\lambda_{E,4}$ in the proof of Proposition~\ref{prop:twisted family}, replacing the gluing maps induced by $\tilde{f}$ with those induced by $\tau\circ\tilde{f}$ does not change the induced gluing maps on $R$ or $\lambda_4$. Consequently, the identity maps on the corresponding line bundles over $[0,1]$ descend to natural isomorphisms
$\psi^R:R_E^+\to R_E^-$, and $\psi^\lambda:\lambda_{E,4}^+\to \lambda_{E,4}^-$.
Moreover, by the same argument used to show that $\phi$ is constant with respect to the product trivializations of $R_E$ and $\lambda_{E,4}$ in the proof of Proposition~\ref{prop:twisted family}, the construction of the 8-divisibility framing is compatible with the action of $\tau$. Hence these isomorphisms make the following diagram commute up to homotopy:
\[
  \begin{CD}
    (R_E^+)^{\otimes2} @>{\phi^+}>> \lambda_{E,4}^+ \\
    @V{(\psi^R)^{\otimes2}}VV @V{\psi^\lambda}VV \\
    (R_E^-)^{\otimes2} @>{\phi^-}>> \lambda_{E,4}^-.
  \end{CD}
\]

We now compare the stable framings $\xi_{\Dsla}^+$ and $\xi_{\Dsla}^-$ using the half-rotation $[0,1]\to\GL_{\C}(V^\pm)$ given by $t\mapsto e^{i\pi t}$, which identifies the two families spin$^c$ structures. By the above commutative diagram, as in the proof of Proposition~\ref{prop:twisted family}, the induced identification of determinant lines differs by the loop $t\mapsto e^{i\pi rt}$, where $r=\ind_{\C}\Dsla_{\fraks}$ is even. This loop has winding number $r/2=\sigma(X)/16$ mod 2.
Hence, if $\sigma(X)$ is divisible by 32, we have $\xi_{\Dsla}^+=\xi_{\Dsla}^-$.
This proves the independence of \eqref{eq: FBF in triviality lemma} of extensions $(\frakt_E, L_E)$.
This completes the proof.
\end{proof}

Let $\TDiff(X,\frakt,L)$ denote the subgroup of $\TDiff(X)$ consisting of diffeomorphisms $f$ of $X$ satisfying $f^\ast\frakt\cong\frakt$ and $f^\ast L\cong L$.
If $X$ satisfies the assumption of Proposition~\ref{prop: 8 framing can detect nontriviality}, for $f \in \TDiff(X,\frakt,L)$, the value \eqref{eq: FBF in triviality lemma} in $\pi_{d(\fraks)+2}^{s}$ is independent of these extensions $(\frakt_E,L_E)$.
Thus 
\[
\FBF_{\frakt,L}^{\divi}(f) := \FBF(X,\fraks,f,\xi_{\Dsla}^{\divi}(\frakt_E,L_E),\xi_{H^+}^{\std})
\]
defines a well-defined map 
\begin{align}
\label{eq: FBF map for the special framing}
\FBF_{\frakt,L}^{\divi}: \TDiff(X,\frakt,L) \to \pi_{d(\fraks)+2}^{s}.
\end{align}

\begin{prop}
\label{prop: FBF hom}
Let $X$ be a closed smooth spin 4-manifold with $b_1(X)=0$ and $\sigma(X)$ divisible by $32$.
Let $\frakt$ be a spin structure on $X$, and let $L\to X$ be a complex line bundle.
Set $\fraks=\frakt\otimes L^4$.
Then the map \eqref{eq: FBF map for the special framing} descends to a well-defined group homomorphism
\[
\FBF_{\frakt,L}^{\divi} : \pi_0(\TDiff(X,\frakt,L)) \to \pi_{d(\fraks)+2}^{s}.
\]
\end{prop}

\begin{proof}
Set $G=\TDiff(X,\frakt,L)$.
We show that the map $\FBF_{\frakt,L}^{\divi}:G\to\pi_{d(\fraks)+2}^{s}$ is additive and descends to $\pi_0(G)$.

Let $f,g\in G$.
Choose lifts $\tilde{f}_{\frakt} : \frakt \to \frakt$, $\tilde{f}_{L} : L \to L$ of $f$ to $\frakt$
and to $L$, respectively.
Similarly, we choose lifts $\tilde{g}_{\frakt} :\frakt \to \frakt$, $\tilde{g}_{L} : L \to L$ of $g$.
By identifying the distinguished fibers of their mapping tori, we obtain a family $E\to B:=S^1\vee S^1$ together with extensions $(\frakt_E,L_E)$.
Since both monodromies act trivially on $H^2(X;\R)$, $H^+(E)$ has the standard framing.
Equip the real families Dirac index with the 8-divisibility framing associated with $(\frakt_E,L_E)$.
By naturality, these framings restrict to the prescribed framings on the two circles.

Put $k=d(\fraks)+1$, and let 
\[
\Phi_E := \widetilde{\FBF}
(E,\fraks_E,
\xi_{\Dsla}^{\divi}(\frakt_E,L_E),
\xi_{H^+}^{\std}) : S^k\wedge B_+\to S^0
\]
denote the finite-dimensional approximation for $E$ with respect to the framings $\xi_{\Dsla}^{\divi}(\frakt_E,L_E)$, $\xi_{H^+}^{\std}$.
Use the stable splittings determined by the basepoint inclusions and the augmentations, as in \eqref{eq:stable spliting}.
Under the splitting
\[
S^k\wedge B_+ \simeq S^k\vee S^{k+1}\vee S^{k+1},
\]
the components of $\Phi_E$ are
$\BF(X,\fraks),\ 
\FBF_{\frakt,L}^{\divi}(f),\ 
\FBF_{\frakt,L}^{\divi}(g)$, respectively.

Let $c:S^1\to B$ be the concatenation of the two loops, with monodromy $fg$.
Then $c^*E$ is the mapping torus of $fg$, and the induced extensions of $\frakt$ and $L$ are given by the compositions $\tilde{f}_{\frakt} \circ \tilde{g}_{\frakt}$ and $\tilde{f}_{L} \circ \tilde{g}_{L}$.
Let $c_+:S^1_+\to B_+$ be the map obtained from $c$ by adjoining disjoint basepoints to its domain and codomain.
By the naturality of the families Bauer--Furuta map and the 8-divisibility framing and the standard framing for $H^+$, the finite-dimensional approximation for $c^*E$ is given by $\Phi_E\circ(\id_{S^k}\wedge c_+)$, of which
the second component along the splitting \eqref{eq:stable spliting} is $\FBF_{\frakt,L}^{\divi}(f)+\FBF_{\frakt,L}^{\divi}(g)$.
By the independence of the extensions, this equals $\FBF_{\frakt,L}^{\divi}(fg)$.
Hence we obtain the desired additivity
\[
\FBF_{\frakt,L}^{\divi}(fg)
=
\FBF_{\frakt,L}^{\divi}(f)+\FBF_{\frakt,L}^{\divi}(g).
\]
This proves the additivity of $\FBF_{\frakt,L}^{\divi} : G \to \pi_{d(\fraks)+2}^{s}$.

Next, we prove $\FBF_{\frakt,L}^{\divi}(f)=0$ when $f \in G$ is smoothly isotopic to the identity.
In this case, we identify $E=E_f$ with $X\times S^1$ and choose $\frakt_E$ and $L_E$ to be the product families $\underline{\frakt}=\frakt\times S^1$ and $\underline{L}=L\times S^1$, respectively.
For these choices, we may take the finite-dimensional approximation $\Psi_E$, given as \eqref{eq: fin dim approx Psi}, for $(X\times S^1,\underline{\frakt}\otimes\underline{L}^4)$ to be $\Psi_X\times\id_{S^1}$.
Since the 8-divisibility framing agrees with the product framing and we have already shown that the invariant is independent of the choices of extensions $(\frakt_E,L_E)$, we obtain $\FBF_{\frakt,L}^{\divi}(f)=0$.

By this triviality for the identity element of $\pi_0(G)$ and the additivity of $\FBF_{\frakt,L}^{\divi} : G \to \pi_{d(\fraks)+2}^{s}$ proven above, $\FBF_{\frakt,L}^{\divi}$ descends to a homomorphism $\FBF_{\frakt,L}^{\divi} : \pi_0(G) \to \pi_{d(\fraks)+2}^{s}$.
This completes the proof.
\end{proof}

If, in addition, $X$ is simply-connected, then every element of $\TDiff(X)$ preserves $\frakt$ and $L$, hence $\TDiff(X)=\TDiff(X,\frakt,L)$.
In this case, we write $\FBF_{\fraks}^{\divi}=\FBF_{\frakt,L}^{\divi}$ for $\fraks=\frakt \otimes L^4$.
Namely, for $f \in \TDiff(X)$, we define
\[
\FBF_{\fraks}^{\divi}(f)
= \FBF(X,\fraks,f,\xi_{\Dsla}^{\divi}(\frakt_E,L_E),\xi_{H^+}^{\std}),
\]
which defines a map
\begin{align}
\label{eq: FBF divi simplyconnected}
\FBF_{\fraks}^{\divi} : \TDiff(X) \to \pi_{d(\fraks)+2}^{s}.
\end{align}

\begin{cor}
\label{cor: FBF hom simply connected}
Let $X$ be a closed simply-connected smooth spin 4-manifold with $\sigma(X)$ divisible by $32$.
Let $\fraks$ be a spin$^c$ structure on $X$ such that $c_1(\fraks)$ is divisible by 8.
Then the map \eqref{eq: FBF divi simplyconnected} descends to a well-defined group homomorphism
\[
\FBF_{\fraks}^{\divi} : \pi_0(\TDiff(X)) \to \pi_{d(\fraks)+2}^{s}.
\]    
\end{cor}

\begin{proof}
This is an obvious consequence of Proposition~\ref{prop: FBF hom}.
\end{proof}

\subsection{Boundary Dehn twists of punctured 4-manifolds}
\label{subsection Boundary Dehn twists of punctured 4-manifolds}

We now turn to Dehn twists along embedded 3-spheres.
Let $X$ be a smooth 4-manifold, and let $Y\subset X$ be a smoothly embedded 3-sphere with trivial normal bundle.
Fix an identification of a tubular neighborhood $\nu(Y)$ with $S^3\times[0,1]$.
Let $\alpha:[0,1]\to\SO(4)$ be a smooth loop representing the nontrivial element of $\pi_1(\SO(4))\cong\Z/2$, such that $\alpha(t)$ is the identity near $t=0,1$.
The {\it Dehn twist} along $Y$ is the diffeomorphism $\delta:X\to X$ defined on $\nu(Y)$ by $\delta(x,t)=(\alpha(t)x,t)$ for $(x,t)\in S^3\times[0,1]$ and extended by the identity outside $\nu(Y)$.

If $X$ is a closed 4-manifold, put $\mathring{X}=X\setminus\mathrm{Int}(D^4)$.
The \emph{boundary Dehn twist} $\delta_{\mathring{X}}\in\Diff_\del(\mathring{X})$ is the Dehn twist along a 3-sphere parallel to $\del\mathring{X}$, supported in a collar neighborhood of the boundary.

As a consequence of the computation for the twisted family, Proposition~\ref{prop:twisted family}, we obtain the following result on the detection of boundary Dehn twists.
Recall that we defined $\overline{\sigma}(X) := \frac{\sigma(X)}{16} \bmod 2$.

\begin{thm}
\label{thm: puncture boundary Dehn twist}
Let $X$ be a simply-connected closed smooth spin 4-manifold.
Suppose that there exists a spin$^c$ structure $\fraks$ on $X$ for which $c_1(\fraks)$ is divisible by 8.
If 
\[
\overline{\sigma}(X) \cdot \eta \cdot \BF(X,\fraks) \neq 0
\]
in $\pi_{d(\fraks)+2}^{s}$, then the boundary Dehn twist $\delta_{\mathring{X}}$ is nontrivial in $\pi_0(\Diff_\del(\mathring{X}))$.    
\end{thm}

\begin{proof}
Let $\frakt$ denote the unique spin structure on $X$.
Then $\fraks$ can be written as $\fraks=\frakt \otimes L^4$ for some complex line bundle $L \to X$.

Let $(X^\tau,\frakt^\tau)$ denote the mapping torus of $\tau$, whose underlying manifold bundle is $X \times S^1$.
Assume that $\delta_{\mathring{X}}$ is trivial in $\pi_0(\Diff_\del(\mathring{X}))$.
Then, by the argument in \cite[Section 2]{KM20Dehn}, there exists a based loop $g_t$ in $\Diff^+(X)$ whose lift to the automorphism group of $(X,\frakt)$ starts at the identity
and ends at $\tau$.
Consequently, the bundle automorphism $F:X\times S^1\to X\times S^1$ defined by $F(x,[t])=(g_t(x),[t])$ induces an isomorphism $(X^\tau,\frakt^\tau)
\cong (X,\frakt)\times S^1$ as families of spin 4-manifolds over $S^1$.

Take $\underline{L}=L\times S^1$ as the extension of $L$ to $X\times S^1$.
Since $H^1(X;\Z)=0$, the restriction map $H^2(X\times S^1;\Z)\to H^2(X;\Z)$ is an isomorphism.
As each $g_t$ is isotopic to the identity, it follows that $F^*\underline{L}\cong\underline{L}$.
Thus $F$ identifies the corresponding pairs of families spin structures and line bundles.

By the naturality of the 8-divisibility framing and the definition of the standard $H^+$-framing, this family isomorphism preserves both framings.
Hence we have
\[
\FBF(X,\fraks,\id_X,
\xi_{\Dsla}^{\divi}(\underline{\frakt}, \underline{L}),
\xi_{H^+}^{\std})
=
\FBF(X,\fraks,\id_X,
\xi_{\Dsla}^{\divi}(\frakt^\tau, \underline{L}),
\xi_{H^+}^{\std}),
\]
where $\underline{\frakt}=\frakt\times S^1$.
By construction,
$\xi_{\Dsla}^{\divi}(\underline{\frakt},\underline{L})$
is the product framing, so the left-hand side is zero.
Thus, we have
\begin{align}
\label{eq: punc boundary FBF zero}
\FBF(X,\fraks,\id_X,
\xi_{\Dsla}^{\divi}(\frakt^\tau,\underline{L}),
\xi_{H^+}^{\std})=0.
\end{align}
On the other hand, by Proposition~\ref{prop:twisted family} and our assumption, we have
\[
\FBF(X,\fraks,\id_X,
\xi_{\Dsla}^{\divi}(\frakt^\tau,\underline{L}),
\xi_{H^+}^{\std})
=
\overline{\sigma}(X)\cdot\eta\cdot\BF(X,\fraks)
\neq 0
\]
in $\pi^s_{d(\fraks)+2}$, which contradicts \eqref{eq: punc boundary FBF zero}.
\end{proof}

\subsection{Gluing}

To handle Dehn twists on connected sums, we need the following gluing formula.
Let $X_1,X_2$ be closed smooth spin 4-manifolds with $b_1(X_i)=0$. 
For $i=1,2$, let $\frakt_i$ be a spin structure on $X_i$ and let $L_i\to X_i$ be a complex line bundle.
Set $\fraks_i=\frakt_i\otimes L_i^4$.
We put $X=X_1\#X_2$, $\fraks=\fraks_1\#\fraks_2$, and $L=L_1\#L_2$.
We orient $H^+(X)$ by the direct sum of the fixed orientations of $H^+(X_1)$ and $H^+(X_2)$.
Let $\delta:X\to X$ denote the Dehn twist along the neck of the connected sum $X=X_1\#X_2$.
Note that $\delta$ preserves $\fraks$.
Let
$\tau_i:\frakt_i\to\frakt_i$ denote the nontrivial deck transformation of $\frakt_i$.
Let $E_{\delta}\to S^1$ be the mapping torus of $\delta$.

Let $\frakt_{E_{\delta}}^{X_1}$ denote the families spin structure on $E_\delta$ with fiber $\frakt=\frakt_1\#\frakt_2$ that corresponds to the twisted family on $X_1$ by the deck transformation.
Namely, the restriction of $\frakt_{E_{\delta}}^{X_1}$ to $\mathring{X}_1$ is the twisted family $\frakt_1^{\tau_1}|_{\mathring{X}_1}$, while the restriction of $\frakt_{E_{\delta}}^{X_1}$ to $\mathring{X}_2$ is the trivial family $\underline{\frakt_2}|_{\mathring{X}_2}$. 
To glue them, let $\alpha_t\in \SO(4)$ ($t\in[0,1]$) be a based loop representing the nontrivial element of
$\pi_1(\SO(4))\cong\Z/2$.
Its spin lift $\tilde{\alpha}_t$, starting at the identity, ends at the nontrivial deck
transformation.
The families spin structure $\frakt_{E_\delta}^{X_1}$ is obtained by gluing the twisted family on $X_1$ and the product family on $X_2$ along the neck
using these spin lifts, as in \cite[Section~5]{KM20Dehn}.
%These two families are glued together along the neck using a path lifting a loop representing the generator of $\pi_1(\SO(3))$.

Let $L_{E_\delta}\to E_\delta$ be a complex line bundle whose restriction to each fiber is isomorphic to $L$. Such a line bundle exists and is unique up to isomorphism, since $H^1(X;\Z)=0$ and $\delta$ acts trivially on $H^2(X;\Z)$, so the Wang exact sequence gives an isomorphism $H^2(E_\delta;\Z)\cong H^2(X;\Z)$.

\begin{prop}
\label{prop: FBF gluing}
In the above setup, we have
\[
\begin{aligned}
&\FBF(X,\fraks,\delta,
\xi_{\Dsla}^{\divi}
(\frakt_{E_{\delta}}^{X_1},
L_{E_{\delta}}),
\xi_{H^+}^{\std})
=
\FBF(X_1,\fraks_1,\id_{X_1},
\xi_{\Dsla}^{\divi}
(\frakt_1^{\tau_1},\underline{L}_1),
\xi_{H^+}^{\std})
\cdot\BF(X_2,\fraks_2)
\end{aligned}
\]
in
$\pi_{d(\fraks)+2}^{s}$, where
$d(\fraks)=d(\fraks_1)+d(\fraks_2)+1$.
\end{prop}

\begin{proof}
Equip $E_1=X_1\times S^1$ and $E_2=X_2\times S^1$ with
$(\frakt_1^{\tau_1},\underline{L}_1)$ and
$(\underline{\frakt_2},\underline{L}_2)$, respectively.
We use Bauer's gluing theorem \cite[Theorem~2.1]{B04}, applied fiberwise as in \cite[Proof of Proposition~5.1]{KM20Dehn}, and check its compatibility with the prescribed framings.

Following \cite[Proof of Proposition~5.1]{KM20Dehn},
consider families $Z,Z'\to S^1$ with fibers
\[
(X_1\#S^4)\sqcup(S^4\#X_2)\sqcup(S^4\#S^4),
\quad
(X_1\#X_2)\sqcup(S^4\#S^4)\sqcup(S^4\#S^4),
\]
respectively.
The $X_i$-pieces of the first two components of $Z$ carry the triples $(E_1,\frakt_1^{\tau_1},\underline{L}_1)$ and  $(E_2,\underline{\frakt_2},\underline{L}_2)$, respectively.
The additional $S^4$-pieces are auxiliary.
The three neck collars are cyclically reglued, as in \cite[Proof of Proposition~5.1]{KM20Dehn}, so that the first component of $Z'$ is $E_\delta$
with the families spin structure $\frakt_{E_\delta}^{X_1}$.

On the collar of the component $X_1\# S^4$, the self-diffeomorphism of $S^1\times[-1,1]\times S^3$ defined by 
$(t,r,y)\mapsto(t,r,\alpha_t y)$ identifies its twisted spin structure with the product spin structure.
%Indeed, the spin lift of $\alpha_t$ over $[0,1]$ ends at the deck transformation $\tau_1$, which agrees with the endpoint identification of the twisted spin family.
%Together with compatible product identifications on the remaining collars, this gives a common product spin model for all three collar families.
On the other two collars, lying in the $S^4\#X_2$ and $S^4\#S^4$ components, we use the product spin identifications. Together with the identification of the first collar described above, these are chosen so that the cyclic permutation of the three collars lifts to an isomorphism of families spin bundles over $S^1$, as required in \cite[Proof of Proposition~5.1]{KM20Dehn}.

%Choose trivializations of $L_i$ over the boundary 3-spheres used in the connected-sum construction, and extend them constantly over the base $S^1$.
%Such trivializations exist and are unique up to homotopy on each $S^3$, since $H^2(S^3;\Z)=0$ and $\pi_3(U(1))=0$.
%The restrictions of the coefficient line bundles to the collar families are also trivial, since $H^2(S^1\times[-1,1]\times S^3;\Z)=0$.
%We choose their trivializations to be compatible with the spin collar identifications above.
%Using these trivializations and identity gluing maps on the coefficient line bundles $L_i$'s, we obtain line bundles $L_Z$ and $L_{Z'}$ on $Z$ and $Z'$, respectively.

Choose trivializations of $L_i$ over the boundary 3-spheres used in the connected-sum construction, and extend them over the collars and constantly over the base $S^1$.
Such trivializations exist and are unique up to homotopy on each $S^3$, since $H^2(S^3;\Z)=0$ and $\pi_3(U(1))=0$.
On the auxiliary $S^4$-pieces, use trivial line bundles with their standard trivializations.
Gluing $L_i$'s using these trivializations, we obtain line bundles $L_Z$ and $L_{Z'}$ on $Z$ and $Z'$, respectively.
By construction, their collar trivializations are compatible with the cyclic regluing.

Choose compatible product metrics and unitary trivializations of $L_Z$ and $L_{Z'}$ on the collars, and take their connections to be trivial there.
Choose the section $v_Z$ of $L_Z^{\oplus4}$ used in Lemma~\ref{lem: 8-divisibility framing} to have the same nonzero constant value, say $(1,0,0,0) \in \C^4$, on all collars, and let $v_{Z'}$ be its regluing.
The prescribed collar section extends without zeros because $\dim Z=5$ and the fibers of the sphere bundle $S(L_Z^{\oplus4})$ are 6-connected.

Recall that, for a fixed $k\geq4$, an isomorphism  
\[
V : L^2_{k}(Z|_{\{\pt\}};S^\pm) \to  L^2_{k}(Z'|_{\{\pt\}};S^\pm)
\]
is introduced in \cite[page 25]{B04}. In compatible collar identifications, these maps multiply the vector of sections on the three collars by $A(\rho)$, where $A:[0,1]\to\SO(3)$ joins the identity to the cyclic permutation matrix and $\rho : S^3 \times [-L,L] \to [0,1]$ for $L>0$ is a smooth collar function, constant in the $S^3$-direction, equal to $0$ and $1$ near the two ends of $[-L,L]$.
Outside the collars, they identify the unchanged pieces.

Write $D_j=\Dsla_{Z,j}$ and $D'_j=\Dsla_{Z',j}$ for
$j=0,\ldots,4$.
By extending Bauer's construction of $V$ to the family setting and to $L^j$-twisted spinor bundles, we obtain isomorphisms of Hilbert bundles over $B$:
\[
\begin{aligned}
V_j^+ &: L^2_{k,B}(Z;S_Z^+\otimes L_Z^j) \to L^2_{k,B}(Z';S_{Z'}^+\otimes L_{Z'}^j),\\
V_j^- &:L^2_{k-1,B}(Z;S_Z^-\otimes L_Z^j) \to  L^2_{k-1,B}(Z';S_{Z'}^-\otimes L_{Z'}^j).
\end{aligned}
\]
Here $L^2_{k,B}(-)$ denotes the fiberwise Sobolev completion of $\Gamma_B(-)$.
Use the same $A$ and $\rho$ for every $j$.

Let $Q_Z^\pm$ and $Q_{Z'}^\pm$ be the maps induced by $Q_v$ on positive and negative spinors.
Using the exterior-power decompositions in the proof of Lemma~\ref{lem: 8-divisibility framing}, form the corresponding direct sums $V_{\even}^\pm$ and $V_{\odd}^\pm$.
Since $A(\rho)$ mixes collar components while $Q_v$ acts on their common exterior-algebra factor, we have
\begin{align}
\label{eq: intertwining}
V_{\odd}^\pm \circ Q_Z^\pm = Q_{Z'}^\pm \circ V_{\even}^\pm.
\end{align}

Set
\[
\widehat{D}'_j=(V_j^-)^{-1}\circ D'_j\circ V_j^+: L^2_{k,B}(Z;S_Z^+\otimes L_Z^j) \to L^2_{k-1,B}(Z;S_Z^-\otimes L_Z^j)
\]
and
\[
D_{j,t} = (1-t)D_j + t\widehat{D}'_j
\]
for $t\in[0,1]$.
The difference $\widehat{D}'_j-D_j$ is supported on the collars and is zeroth-order in the collar identifications, hence is a compact operator between these Sobolev spaces.
Thus $D_{j,t}$ is a homotopy of families Fredholm operators.
For $j=4$, it gives the index identification in Bauer's gluing argument: its path and the spinor linear path in that argument lie in the same affine space of compact perturbations of $D_4$, and hence are homotopic relative to their endpoints \cite[Section~3]{B04}.

Put
\[
D_t^{\even} = D_{0,t}\oplus D_{2,t}^{\oplus6}\oplus D_{4,t},\ 
D_t^{\odd} = D_{1,t}^{\oplus4}\oplus D_{3,t}^{\oplus4}.
\]
Consider a 2-parameter family
\[
H_{s,t} = (1-s)Q_Z^- \circ D_t^{\even} \circ (Q_Z^+)^{-1}
+sD_t^{\odd}
\]
parameterized by $(s,t)\in[0,1]^2$.
Each $H_{s,t}$ is Fredholm, since it is a compact perturbation
of $D_{t=0}^{\odd}$.
By \eqref{eq: intertwining}, the restrictions of $\{H_{s,t}\}_{(s,t) \in [0,1]^2}$ to $t=0,1$ are the homotopies \eqref{eq: Fred homotopy} used in Lemma~\ref{lem: 8-divisibility framing} for $Z$ and $Z'$, after transfer to $Z$ by conjugation by $V_{\odd}^\pm$.
The restrictions of $H_{s,t}$ to $s=0$ and $s=1$ are $Q_Z^-\circ D_t^{\even} \circ (Q_Z^+)^{-1}$ and $D_t^{\odd}$, respectively. Thus, they correspond to the homotopies induced by the gluing identifications for the even and odd exterior-power twists, with the former transported to the odd Sobolev bundles via $Q_Z^\pm$.
%Since $H_{s,t}$ defines a family of Fredholm operators over $S^1\times[0,1]^2$, the two paths along the boundary of $[0,1]^2$ from $(0,0)$ to $(1,1)$ induce homotopic isomorphisms between determinant line bundles.

For $j=0,\ldots,4$, let $T_j:\lambda_{Z,j}\to\lambda_{Z',j}$ denote the isomorphism induced by the homotopy $D_{j,t}$ and the identifications $V_j^\pm$.
Set
\[
T_R= T_3^{\otimes2}\otimes T_2^{\otimes(-3)}\otimes T_1^{\otimes2} : R_Z\to R_{Z'}.
\]
Moreover, $V_0^\pm$ and $D_{0,t}$ are quaternionic linear, since the collar identifications preserve the spin structures and $A(\rho)$ has real entries.
Thus $T_0$ preserves the canonical quaternionic trivializations of $\lambda_{Z,0}$ and $\lambda_{Z',0}$.
The 2-parameter family $\{H_{s,t}\}$ of (families over $B=S^1$ of) Fredholm operators therefore implies
\begin{align}
\label{eq: compatibility gluing}
T_4\circ\phi_Z \simeq \phi_{Z'}\circ T_R^{\otimes2}.    
\end{align}

On the components of $Z$ and $Z'$ with fibers diffeomorphic to $S^4$, choose positive scalar curvature metrics, trivial connections on the restrictions of $L_Z$ and $L_{Z'}$, and constant sections $v$. 
Their twisted Dirac operators are invertible, so the determinant line bundles $\lambda_{E,j}$ are canonically trivial.
The square root data $(R_E,\phi_E)$ constructed in Lemma~\ref{lem: 8-divisibility framing} therefore agree, up to homotopy, with the trivial square root data.
Moreover, as in \cite[Proof of Proposition~5.1]{KM20Dehn}, the finite-dimensional approximation of the families Seiberg--Witten map of each such $S^4$-family, for either choice of family spin structure, represents the multiplicative unit.

Combining the compatibility \eqref{eq: compatibility gluing} with multiplicativity under disjoint union, we obtain an isomorphism of square root data
\[
(R_{E_\delta},\phi_{E_\delta})
\cong
(R_{E_1}\otimes R_{E_2}, \phi_{E_1}\otimes\phi_{E_2})
\]
up to homotopy, compatible with the identifications of the determinant line bundles obtained from $T_j$.
Consequently, the corresponding stable identification of the real families Dirac indices preserves the spin structures defined by Lemma~\ref{lem:Sasahira} determined by the square root data $(R_{E_\delta},\phi_{E_\delta})$ and $(R_{E_1}\otimes R_{E_2}, \phi_{E_1}\otimes\phi_{E_2})$
and hence, by Lemma~\ref{lem:spin induces framing},
the 8-divisibility framings up to homotopy.
%Consequently, the gluing identification preserves the induced spin structures on the real families Dirac indices and, by Lemma~\ref{lem:spin induces framing}, their 8-divisibility framings.
On the other hand, the standard framing of $H^+(E_\delta)$ agrees, up to homotopy, with the direct sum of the standard framings of $H^+(E_1)$ and $H^+(E_2)$, under the obvious isomorphism $H^+(E_\delta) \cong H^+(E_1)\oplus H^+(E_2)$.

%On the components of $Z$ and $Z'$ with fibers diffeomorphic to $S^4$, choose positive scalar curvature metrics, trivial connections on the restrictions of $L_Z$ and $L_{Z'}$, and constant sections $v$.
%Their Dirac operators are invertible, so all the determinant line bundles $\lambda_{E,j}$ are canonically trivial. Since the sections $v$ are constant and the connections are trivial, the isomorphism $\phi_E$ constructed in Lemma~\ref{lem: 8-divisibility framing} agrees, up to homotopy, with the canonical isomorphism between the trivial line bundles $R_E^{\otimes 2}$ and $\lambda_{E,4}$.
%By \cite[Proof of Proposition~5.1]{KM20Dehn}, the families Bauer--Furuta invariant of each such $S^4$-family, for either choice of family spin structure, is represented by the identity map between zero-dimensional vector bundles over $S^1$, and hence is the multiplicative unit. 

The gluing theorem therefore gives the fiberwise smash product formula with these framings.
Since $E_2$ has product data and product framing, its finite-dimensional approximation can be chosen constant over $S^1$.
Taking the second component under \eqref{eq:stable spliting} gives the asserted formula.
\end{proof}

\subsection{Dehn twists on connected sums}

Now, we are ready to detect Dehn twists on connected sums.
If $X=X_1\#X_2$ is a connected sum of closed smooth 4-manifolds, the {\it Dehn twist along the neck} is the Dehn twist along the separating 3-sphere arising from the connected-sum construction.

\begin{thm}
\label{thm: Dehn twists on connected sums}
Let $X_1,X_2$ be closed smooth spin 4-manifolds with $b_1(X_i)=0$ and $\sigma(X_i)\equiv 16 \bmod 32$. 
For $i=1,2$, let $\frakt_i$ be a spin structure on $X_i$ and let $L_i\to X_i$ be a complex line bundle.
Set $\fraks_i=\frakt_i\otimes L_i^4$.
Assume that, for $i=1,2$, $d(\fraks_i)=0$ and  $\SW(X_i,\fraks_i)$ is odd.
Then the Dehn twist $\delta$ on $X_1\#X_2$ along the neck is not smoothly isotopic to the identity.
\end{thm}

\begin{proof}
First, note that the divisibility of $c_1(\fraks_i)$ by 8 and $d(\fraks_i)=0$ implies that $b^+(X_i)>1$, as described in Lemma~\ref{lem: BF equals eta iff SW odd}, so we can define $\SW(X_i,\fraks_i)$.
By Lemma~\ref{lem: BF equals eta iff SW odd}, we have $\BF(X_i,\fraks_i) = \eta$.
Since $\sigma(X_1) \equiv 16$ mod 32, it follows from Proposition~\ref{prop:twisted family} that
\[
\FBF(X_1,\fraks_1,\id_{X_1},\xi_{\Dsla}^{\divi}(\frakt_1^{\tau_1}, \underline{L_1}),\xi_{H^+}^{\std})
= \eta^2
\]
in $\pi^{s}_{2}$.
This combined with the gluing result, Proposition~\ref{prop: FBF gluing}, implies that 
\[
\FBF(X_1\#X_2,\fraks_1\#\fraks_2,\delta, \xi_{\Dsla}^{\divi}(\frakt_{E_{\delta}}^{X_1}, L_{E_{\delta}}), \xi_{H^+}^{\std})
= \eta^3 \neq 0
\]
in $\pi^{s}_{3}$,
where we used the notation in Proposition~\ref{prop: FBF gluing}.
Since $\sigma(X_1\#X_2)$ is divisible by 32, we may apply Proposition~\ref{prop: FBF hom} to $X_1\#X_2$ and conclude that $\delta$ is not smoothly isotopic to the identity.
\end{proof}

\begin{proof}[Proof of Theorem~\ref{thm: eg Dehn twists on connected sums}]
The diffeomorphism $\delta$ is not smoothly isotopy to the identity by Theorem~\ref{thm: Dehn twists on connected sums} and Remark~\ref{rem: simplyconn 8 divi}, but is topologically isotopic to the identity by a result of Quinn~\cite{Q86}, together with the correction by Gabai--Gay--Hartman--Krushkal--Powell~\cite{gabai2023pseudoisotopies}.
\end{proof}

\subsection{Detecting Dehn twists from their square roots}

Next, we detect Dehn twists from the actions of their square roots on homology.
This is typically useful to detect Dehn twists on irreducible 4-manifolds.
The argument is inspired by that of Miyazawa~\cite{miyazawaDehn}.
Miyazawa considered spin 4-manifolds with boundary and used quaternionic framings.
By using the 8-divisibility framing in place of the quaternionic framing, we extend his argument to non-spin spin$^c$ structures and to closed 4-manifolds.

Again, note that the existence of a spin$^c$ structure $\fraks$ such that $c_1(\fraks)$ is divisible by 8 and $d(\fraks)=0$ implies that $b^+(X)>1$.
We denote by $H^+(X)$ a choice of maximal-dimensional positive-definite subspace of $H^2(X;\R)$.
As mentioned, the space of such choices is contractible.
Suppose that there is an automorphism $A$ of the intersection form of $X$ such that $A^2=\id$.
Then there is a choice of $H^+(X)$ such that $H^+(X)$ is invariant under $A$ setwise.
Indeed, the eigenspaces of $A$ give an orthogonal (with respect to the intersection form) decomposition
\[
H^2(X;\R)=H(A;+1)\oplus H(A;-1),
\]
for the eigenvalues $1, -1$, and the restriction of the intersection form to each eigenspace is nondegenerate.
Choose a maximal-dimensional positive-definite subspace
$P_\pm\subset H(A;\pm1)$.
Then $H^+(X):=P_+\oplus P_-$ is a maximal-dimensional positive-definite subspace of $H^2(X;\R)$ that is setwise invariant under $A$.

\begin{thm}
\label{thm: square root obstruction}
Let $X$ be a closed smooth spin 4-manifold with $b_1(X)=0$ and $\sigma(X)$ divisible by 32.
Let $\frakt$ be a spin structure on $X$, and let $L\to X$ be a complex line bundle.
Set $\fraks=\frakt\otimes L^4$.
Assume that $d(\fraks)=0$ and $\SW(X,\fraks)$ is odd.
Let $f:X\to X$ be an orientation-preserving diffeomorphism satisfying $f^\ast\frakt\cong\frakt$, $f^\ast L\cong L$, and $f^2\in\TDiff(X)$.
Let $H^+(X)$ be a maximal positive-definite subspace of $H^2(X;\R)$ invariant under $f^\ast$,
and put
\[
m=\dim H^+(X)^{-f^\ast}.
\]
Then $m$ is even, and
\[
\FBF_{\frakt,L}^{\divi}([f^2])
=\frac{m}{2}\eta^2
\]
in $\pi_2^s\cong\Z/2$.

In particular, if $m\equiv2\bmod4$, then $\FBF_{\frakt,L}^{\divi}([f^2])=\eta^2\neq0$, so $f^2$ is not smoothly isotopic to the identity.
\end{thm}

\begin{proof}
Since $f^\ast\fraks\cong\fraks$, we have $\SW(X,\fraks)=(-1)^m\SW(X,\fraks)$.
As $\SW(X,\fraks)$ is odd, $m$ is even.
Put $n=m/2$.

Let $E=E_f\to S^1$ be the mapping torus of $f$, and choose extensions $(\frakt_E,L_E)$ of $(\frakt,L)$ to $E$.
Equip the families Dirac index with the 8-divisibility framing $\xi_{\Dsla}^{\divi}(\frakt_E,L_E)$.
Fix a linear isomorphism $H^+(X)^{-f^\ast} \cong \R^{2n}$, and consider the half-rotation on each copy of $\R^2$.
Namely, consider the path $[0,1] \to \GL_{\R}(2)^{\times n} \subset \GL_{\R}(2n)$ given by $t \mapsto (e^{it\pi},\cdots,e^{it\pi})$.
This gives a framing of $H^+(E)$, which we denote by $\xi_{H^+}^{\pi}$.
%As in the construction in the original proof, let $\xi_{H^+}^{\pi}$ denote the framing obtained by a half-rotation on each of the $n$ copies of $\R^2$ in $H^+(X)^{-f^\ast}$.

Let $p:S^1\to S^1$ be a map of degree 2.
By Lemma~\ref{lem: pullback}, we have
\begin{align}
\label{eq: FBF zero twice}
\FBF(X,\fraks,f^2, p^\ast\xi_{\Dsla}^{\divi}(\frakt_E,L_E), p^\ast\xi_{H^+}^{\pi})
= 2\FBF(X,\fraks,f, \xi_{\Dsla}^{\divi}(\frakt_E,L_E),\xi_{H^+}^{\pi}) =0
\end{align}
in $\pi_2^s\cong\Z/2$.
By the naturality of the 8-divisibility framing
under pullback and
Proposition~\ref{prop: 8 framing can detect nontriviality}, we have
\begin{align}
\label{eq: FBF translate}
\FBF(X,\fraks,f^2, p^\ast\xi_{\Dsla}^{\divi}(\frakt_E,L_E),\xi_{H^+}^{\std})
=\FBF_{\frakt,L}^{\divi}([f^2]).
\end{align}
On the other hand, $p^\ast\xi_{H^+}^{\pi}$ differs from $\xi_{H^+}^{\std}$ by $n\bmod2$ in $\pi_1(\SO)\cong\Z/2$.
Thus, by Lemma~\ref{lem:framing_change}, Lemma~\ref{lem: BF equals eta iff SW odd}, \eqref{eq: FBF zero twice} and \eqref{eq: FBF translate}, we obtain
\[
0 =\FBF_{\frakt,L}^{\divi}([f^2]) +n\eta\cdot\BF(X,\fraks)
=\FBF_{\frakt,L}^{\divi}([f^2]) +n\eta^2.
\]
Since $\pi_2^s\cong\Z/2$, this gives
\[
\FBF_{\frakt,L}^{\divi}([f^2])
=n\eta^2=\frac{m}{2}\eta^2.
\]
The final assertion follows from Proposition~\ref{prop: FBF hom}.
\end{proof}

\section{Examples}
\label{section Examples}

We now present concrete examples of Dehn twists that can be detected using the results established in the preceding sections.

%In this section, we focus on closed simply-connected smooth 4-manifolds $X$ and punctured 4-manifolds $\mathring{X}$.
%In these settings, a homeomorphism acting trivially on homology is topologically isotopic to the identity, relative to the boundary in the punctured case. This follows from Quinn~\cite{Q86}, together with the correction by Gabai--Gay--Hartman--Krushkal--Powell~\cite{gabai2023pseudoisotopies}, in the closed case, and from Orson--Powell~\cite{OrsonPowell2022} in the punctured case.
%Consequently, in either setting, a homologically trivial diffeomorphism that is not smoothly isotopic to the identity is exotic.

\subsection{The Seiberg--Witten invariant of elliptic surfaces}

Our main examples are elliptic surfaces and their connected sums.
We recall some basic facts about elliptic surfaces, following~\cite[Section~8.3]{GS99}.
Let $E(n)$ denote the simply-connected relatively minimal elliptic surface without multiple fibers with holomorphic Euler characteristic $n$.
For relatively prime positive integers $p,q$, let $E(n)_{p,q}$ denote the elliptic surface obtained from $E(n)$ by performing logarithmic transformations of multiplicities $p$ and $q$ along two regular fibers.
In particular, $E(n)_{1,1}=E(n)$.
When $n\geq2$ is even and $p,q$ are odd, $X=E(n)_{p,q}$ is simply-connected and spin, and
\[
\chi(X)=12n,\quad
\sigma(X)=-8n,\quad
b^+(X)=2n-1.
\]

We will repeatedly use the following computation.

\begin{lem}
\label{lem: elliptic 8-divisible basic classes}
Let $n\geq 2$ be even, and let $p,q$ be relatively prime positive odd integers.
Set $X=E(n)_{p,q}$.
Then $X$ admits a spin$^c$ structure $\fraks$ satisfying
\[
c_1(\fraks) \in 8H^2(X;\Z),\quad
d(\fraks)=0,\quad
\SW(X,\fraks) \equiv 1 \bmod 2
\]
if and only if at least one of the following conditions holds:
\begin{enumerate}
\item $p=q=1$ and $n \equiv 2\bmod 8$;
\item $\{p,q\} = \{1,3\}$ and $n \not\equiv 6 \bmod 8$;
\item $\max\{p,q\} \geq 5$.
\end{enumerate}
\end{lem}

\begin{proof}
Let $u\in H^2(X;\Z)$ be the primitive class such that the Poincar\'e dual of a regular fiber is $pqu$.
Set $N = n-2$ and $C = npq-p-q$.
By \cite[Lecture~2, Section~9]{FS6lec}, the Seiberg--Witten polynomial of $X$ is, up to an overall sign,
\[
\begin{aligned}
\SW_X(t)
=\frac{(t^{pq} - t^{-pq})^n}{(t^p - t^{-p})(t^q - t^{-q})}
=\sum_{j=0}^{N}
  \sum_{a=0}^{q-1}
  \sum_{b=0}^{p-1}
  (-1)^j\binom{N}{j}
  t^{C-2(jpq + pa + qb)},
\end{aligned}
\]
where $t^v$ corresponds to the class $vu$.

The exponents in this sum are pairwise distinct.
Indeed, if
\[
jpq + pa + qb = j'pq + pa' + qb',
\]
reduction modulo $q$ and modulo $p$, together with $\gcd(p,q)=1$ and the ranges of $a,b$, gives $a=a'$ and $b=b'$.
It then follows that $j=j'$.
Thus the coefficient of $t^{C-2(jpq + pa + qb)}$ is $\pm(-1)^j\binom{N}{j}$.

Since $u^2=0$ and $2\chi(X) + 3\sigma(X) = 0$, every basic class has formal dimension zero.
Moreover, since $u$ is primitive, the class
\[
(C-2(jpq + pa + qb))u
\]
is divisible by $8$ if and only if
\[
C-2(jpq + pa + qb) \equiv 0 \bmod8.
\]

Suppose first that $\max\{p,q\}\geq 5$.
Without loss of generality, we may assume $q\geq 5$.
Since $C$ is even and $p$ is odd, there exists $a \in \{0,1,2,3\}$ such that $pa \equiv \frac{C}{2} \bmod4$.
Taking $j=b=0$, we obtain a basic class $(C-2pa)u$ divisible by $8$, with Seiberg--Witten invariant $\pm1$.

It remains to consider $(p,q)=(1,1)$ and $\{p,q\} = \{1,3\}$.
In either case, we shall use the following elementary observation repeatedly: for $e \geq 0$, if $2^e\mid N$ and $\binom{N}{j} \equiv 1 \bmod2$,
then $2^e\mid j$.
Indeed, this follows by comparing coefficients in $(1+t)^N = (1+t^{2^e})^{N/2^e}$ in $\mathbb{F}_2[t]$.

If $p=q=1$, then $a=b=0$ and $C=N$.
Suppose that $N-2j \equiv 0 \bmod8$ and $\binom{N}{j} \equiv 1 \bmod2$.
Since $N$ is even, the observation above gives $2\mid j$, and hence $4 \mid N$.
Applying the observation again gives $4 \mid j$, and therefore $8 \mid N$.
Conversely, if $8 \mid N$, taking $j=0$ gives a basic class divisible by $8$ with Seiberg--Witten invariant $\pm1$.
Thus this case occurs precisely when $n \equiv2 \bmod8$.

Finally, suppose $\{p,q\} = \{1,3\}$.
Without loss of generality, we may assume $(p,q) = (1,3)$.
Then $C = 3n-4$, $0\leq a\leq2$, and $b=0$.
If $n \equiv 0,2,4 \bmod 8$, respectively, take $j=0$ and $a = 2,1,0$.
In each case, $C-2a\equiv 0 \bmod8$, and the corresponding Seiberg--Witten invariant is $\pm1$.

If $n\equiv6 \bmod8$, then $4\mid N$.
For the coefficient indexed by $j$ to be odd, the observation above forces $4\mid j$.
Consequently,
\[
C-2(3j + a) \equiv 6-2a \bmod8.
\]
For $a=0,1,2$, these residues are $6,4,2$, respectively, so none is zero.
This proves the necessity in this case.
\end{proof}

\begin{rem}
As is clear from the proof of Lemma~\ref{lem: elliptic 8-divisible basic classes}, whenever one of conditions (1)--(3) holds, $\fraks$ can be chosen so that $\SW(X,\fraks)=\pm1$.
\end{rem}

\subsection{Examples of Dehn twists on connected sums}

\begin{thm}
\label{thm: elliptic surface sum}
Let $k,l\geq0$, and let $p,q,r,s$ be positive odd integers such that
$\gcd(p,q)=\gcd(r,s)=1$.
Suppose that
\[
\left\{
\begin{aligned}
&\text{either } k \text{ is even or }
\max\{p,q\}\geq5, \text{and}\\
&\text{either } l \text{ is even or }
\max\{r,s\}\geq5.
\end{aligned}
\right.
\]
Then the Dehn twist $\delta$ on
\[
E(4k+2)_{p,q}\#E(4l+2)_{r,s}
\]
along the neck is an exotic diffeomorphism.
\end{thm}

\begin{proof}
Let $X=E(4k+2)_{p,q}$.
Recall that $X$ is simply-connected and spin, and
\[
\sigma(X)=-8(4k+2)\equiv 16\bmod{32}.
\]
If $k$ is even, then $4k+2\equiv2 \bmod8$, so the conditions of Lemma~\ref{lem: elliptic 8-divisible basic classes} are satisfied for any relatively prime positive odd integers $p,q$.
If $k$ is odd, the same conclusion follows from the assumption $\max\{p,q\}\geq5$.

Thus, by Lemma~\ref{lem: elliptic 8-divisible basic classes}, $X$ admits a spin$^c$ structure $\fraks$ such that $c_1(\fraks)\in8H^2(X;\Z)$, $d(\fraks)=0$, and $\SW(X,\fraks)\equiv1 \bmod2$.
The same argument applies to $E(4l+2)_{r,s}$.
The assertion therefore follows from Theorem~\ref{thm: eg Dehn twists on connected sums}.
\end{proof}

\begin{rem}
Although we have shown that the Dehn twist $\delta$ on $X=E(8k+2)\#E(8l+2)$ is nontrivial in $\pi_0(\Diff(X))$, it follows from Lindblad's result \cite[Section~3]{lindblad2026boundarydehntwistscommutators} that the image of $\delta$ in the abelianization $\pi_0(\Diff(X))^{\ab}
\cong H_1(\BDiff(X);\Z)$ is trivial.
\end{rem}

%\begin{cor}
%\label{thm: elliptic surface sum}
%Let $k,l\geq 0$, and let $p,q,r,s$ be positive odd
%integers such that $\gcd(p,q)=\gcd(r,s)=1$.
%Then the Dehn twist $\delta$ on
%\[
%E(8k+2)_{p,q}\#E(8l+2)_{r,s}
%\]
%along the neck is an exotic diffeomorphism.
%\end{cor}

\subsection{Examples of boundary Dehn twists on punctured 4-manifolds}

First, our detection result for Dehn twists along connected-sum necks immediately yields:

\begin{cor}
\label{cor: eg puncture boundary Dehn twist from connected sum}
Let $X$ be a closed simply-connected smooth spin 4-manifold with $\sigma(X)\equiv 16 \bmod 32$. 
Let $\fraks$ be a spin$^c$ structure on $X$ such that $c_1(\fraks)$ is divisible by 8.
Assume that $d(\fraks)=0$ and  $\SW(X,\fraks)$ is odd.
Then the boundary Dehn twist $\delta_{\mathring{X}}$ is relatively exotic.
\end{cor}

\begin{proof}
Apply Theorem~\ref{thm: eg Dehn twists on connected sums} with $X_1=X_2=X$.
This shows that the Dehn twist along the neck of $X\#X$ is not smoothly isotopic to the identity.
Since this Dehn twist is obtained by extending the boundary Dehn twist $\delta_{\mathring X}$ by the identity over the other summand, $\delta_{\mathring X}$ cannot be smoothly isotopic to the identity relative to the boundary.
The topological triviality follows from Orson--Powell~\cite[Theorem E]{OrsonPowell2022}.
\end{proof}

\begin{ex}
\label{ex: boundary Dehn E(8k+2)}
Let $k\geq0$, and let $p,q$ be positive odd integers such that $\gcd(p,q)=1$.
Suppose that either $k$ is even or $\max\{p,q\}\geq5$.
The boundary Dehn twist on the punctured elliptic surface $(E(4k+2)_{p,q})^{\circ}$ is relatively exotic.
In particular, the boundary Dehn twist $\delta$ on $(E(8k+2))^\circ$ is relatively exotic for all integers $k\geq0$.
This extends the result for $K3^\circ=E(2)^\circ$ \cite{BK19}, \cite[Proposition~1.2]{KM20Dehn} and also strengthens a result of Baraglia and the author~\cite[Theorem~1.5]{baraglia2024irreducible4manifoldsadmitexotic}.
\end{ex}

\begin{proof}[Proof of Theorem~\ref{thm: eg puncture boundary Dehn twist}]
The diffeomorphism $\delta_{\mathring{X}}$ is not smoothly isotopic to the identity  by Theorem~\ref{thm: puncture boundary Dehn twist} and Remark~\ref{rem: simplyconn 8 divi}.
The topological triviality follows from \cite{OrsonPowell2022}.
\end{proof}

%\begin{rem}
%Corollary~\ref{cor: eg puncture boundary Dehn twist from connected sum} also directly follows from a variant of Theorem~\ref{thm: puncture boundary Dehn twist} by using a spin$^c$ structure $\fraks$ with $d(\fraks)=0$ and $\BF(X,\fraks)=\eta$, instead of $d(\fraks)=1$ and $\BF(X,\fraks)=\eta^2$.  
%\end{rem}
\begin{rem}
Corollary~\ref{cor: eg puncture boundary Dehn twist from connected sum} also follows directly from Theorem~\ref{thm: puncture boundary Dehn twist}, without appealing to Theorem~\ref{thm: eg Dehn twists on connected sums}, whose proof relies on a gluing argument (Proposition~\ref{prop: FBF gluing}).
Indeed, by Lemma~\ref{lem: BF equals eta iff SW odd}, we have $\BF(X,\fraks)=\eta$, and hence
\[
\overline{\sigma}(X)\cdot\eta\cdot\BF(X,\fraks)
=\eta^2\neq 0.
\]
Thus, the hypothesis of Theorem~\ref{thm: puncture boundary Dehn twist} is satisfied.
\end{rem}

\begin{cor}
\label{cor: elliptic boundary twist families}
Let $n,m\geq2$ be even integers satisfying $n+m\equiv2\bmod4$.
Let $p,q,r,s$ be positive odd integers with $\gcd(p,q)=\gcd(r,s)=1$.
Assume that each of the triples $(n,p,q)$ and $(m,r,s)$ satisfies one of the conditions in Lemma~\ref{lem: elliptic 8-divisible basic classes}.
Then the boundary Dehn twist on
\[
\left(E(n)_{p,q}\#E(m)_{r,s}\right)^\circ
\]
is relatively exotic.
\end{cor}

\begin{proof}
Put $X_1=E(n)_{p,q}$ and $X_2=E(m)_{r,s}$.
Both $X_i$ are simply-connected and spin.
By Lemmas~\ref{lem: elliptic 8-divisible basic classes} and~\ref{lem: BF equals eta iff SW odd}, there are spin$^c$ structures $\fraks_i$ on $X_i$ such that $c_1(\fraks_i)\in8H^2(X_i;\Z)$, $d(\fraks_i)=0$, and $\BF(X_i,\fraks_i)=\eta$.

Set $X=X_1\#X_2$ and $\fraks=\fraks_1\#\fraks_2$.
Then $d(\fraks)=1$, and Bauer's gluing formula \cite[Theorem~1.1]{B04} gives $\BF(X,\fraks)=\eta^2$.
Moreover, $X$ is simply-connected and spin, $c_1(\fraks)$ is divisible by $8$, and $\sigma(X) = -8(n+m) \equiv 16 \bmod32$.
The assertion follows from Theorem~\ref{thm: eg puncture boundary Dehn twist}.
\end{proof}

\begin{thm}
\label{thm: elliptic boundary twist explicit}
Let $k\geq0$ and $l\geq1$, and let $p,q,r,s$ be positive odd integers with
$\gcd(p,q)=\gcd(r,s)=1$.
Assume that $(r,s)\neq(1,1)$ and that either $k$ is even or $\max\{p,q\}\geq5$.
Then the boundary Dehn twist on
\[
\left(E(4k+2)_{p,q}\#E(4l)_{r,s}\right)^\circ
\]
is relatively exotic.
\end{thm}

\begin{proof}
For $n=4k+2$, the conditions in Lemma~\ref{lem: elliptic 8-divisible basic classes} hold precisely when either $k$ is even or $\max\{p,q\}\geq5$.
For $m=4l$, they hold precisely when $(r,s)\neq(1,1)$.
Since $(4k+2)+4l \equiv 2 \bmod4$, the assertion follows from Corollary~\ref{cor: elliptic boundary twist families}.
\end{proof}

\subsection{Seifert fibered Dehn twists on elliptic surfaces}
\label{subsection Seifert fibered Dehn twists on elliptic surfaces}

Now, we detect Dehn twists on elliptic surfaces along Seifert fibered 3-manifolds.
Let $Y$ be a closed Seifert fibered 3-manifold.
The Seifert circle action determines a loop $\rho_t\in\Diff(Y)$, $t\in[0,1]$, based at the identity, which we take to be constant
near $t=0,1$. For an embedding $Y\times[0,1]\subset X$, the {\it Dehn twist along $Y$} is the diffeomorphism $\delta:X\to X$ defined by $\delta(y,t)=(\rho_t(y),t)$ on $Y\times[0,1]$, and by the identity elsewhere.
If $W$ is a compact smooth 4-manifold with $\del W=Y$, the \emph{boundary Dehn twist} $\delta_W\in\Diff_\del(W)$ is defined by the same construction in a collar neighborhood $\nu(\del W)\cong Y\times[0,1]$.

By \cite[Section~8.3]{GS99}, for every $n\geq1$, the elliptic surface $E(n)$ admits a decomposition
\[
E(n)=M(2,3,6n-1)\cup N(n),
\]
where $N(n)$ is the Gompf nucleus. The logarithmic transformations producing $E(n)_{p,q}$ can be performed inside the nucleus, leaving $M(2,3,6n-1)$ unchanged.

For $1\leq i\leq\lfloor n/2\rfloor$, we have $12i-5< 12i-1\leq6n-1$.
By \cite[Exercise~7.3.18]{GS99}, the Milnor fibers
\[
M(2,3,12i-5),\ M(2,3,12i-1)
\]
admit smooth embeddings into $M(2,3,6n-1)$, and hence into $E(n)_{p,q}$, whose images lie outside the nucleus. Moreover, a regular elliptic fiber can be chosen inside the nucleus, disjoint from these embedded Milnor fibers.

\begin{proof}[Proof of Theorem~\ref{thm: Dehn twists on elliptic surfaces}]
Put $X=E(4k)_{p,q}$, and let $M=M(2,3,r)\subset X$, where $r=12i-5$ or $r=12i-1$.

By Lemma~\ref{lem: elliptic 8-divisible basic classes}, there exists a spin$^c$ structure $\fraks=\frakt\otimes L^4$ on $X$ such that $d(\fraks)=0$ and $\SW(X,\fraks)$ is odd.
Moreover, $\sigma(X)=-32k$.
Let $u\in H^2(X;\Z)$ be the primitive class used in the proof of Lemma~\ref{lem: elliptic 8-divisible basic classes}.
Then $c_1(L)$ is an integral multiple of $u$.
As explained above, the standard embedding
$M(2,3,24k-1)\subset X$ is disjoint from the Gompf nucleus, where a regular elliptic fiber lies.
Since $r\leq24k-1$, the standard embedding of $M$ can be taken inside $M(2,3,24k-1)$.
Thus $pqu|_M=0$.
Since $H^2(M;\Z)$ is torsion-free, we obtain $u|_M=0$, and hence $L|_M$ is trivial.

Now, $\gcd(3,r)=1$ and
\[
\frac{(3-1)(r-1)}4=\frac{r-1}{2}
\]
is odd in both cases, $r=12i-5$ and $r=12i-1$.
As described in \cite[Section~5.1 and Proof of Theorem~1.8]{konno2023exotic} and \cite[Lemma~4.1]{miyazawaDehn}, there exists $f\in\Diff_\del(M)$ such that
$f^2$ is the boundary Dehn twist on $M$ and $f^\ast=-\id$ on $H^2(M;\R)$.
Furthermore, we have $b^+(M)\equiv 2 \bmod 4$ (see, for example, \cite[Lemma 4.1]{miyazawaDehn}).

Extend $f$ to $f:X \to X$ by the identity outside $M$.
Then $f^2=\delta$ and $\dim H^+(X)^{-f^\ast}
=b^+(M)\equiv 2 \bmod 4$.
The diffeomorphism $f : X \to X$ preserves $\frakt$, since $X$ is simply-connected and spin, and preserves $L$, since $L|_M$ is trivial.
If $\delta$ were smoothly isotopic to the identity, this contradicts Theorem~\ref{thm: square root obstruction}.
Therefore, $\delta$ is not smoothly isotopic to the identity.

Finally, by \cite{Q86,gabai2023pseudoisotopies}, $\delta$ is topologically isotopic to
the identity.
This completes the proof.
\end{proof}

\begin{rem}
As is clear from the proof of Theorem~\ref{thm: Dehn twists on elliptic surfaces}, our method applies to closed 4-manifolds that are not necessarily elliptic surfaces, provided that the required divisibility conditions and the condition on the Seiberg--Witten invariant are satisfied.
For example, as discussed in~\cite[Subsection~5.2]{konno2025constraintslefschetzfibrationsfourdimensional}, we can consider certain 4-manifolds obtained by knot surgery.
\end{rem}

\begin{rem}
Just as we detected boundary Dehn twists on Gompf nuclei in Corollary~\ref{cor: nuclei}, we can detect boundary Dehn twists on any symplectic cap of $\Sigma(2,3,12i-5)$ or $\Sigma(2,3,12i-1)$, provided that the required divisibility conditions are satisfied.
\end{rem}

\begin{rem}
In~\cite{konno2025constraintslefschetzfibrationsfourdimensional}, the authors gave a negative answer to a question posed by Donaldson concerning the symplectic Torelli group.
Theorem~\ref{thm: Dehn twists on elliptic surfaces} provides further examples yielding a negative answer to this question.
See~\cite{konno2025constraintslefschetzfibrationsfourdimensional} for details.
\end{rem}

\section{Questions}

We close the paper by posing further questions that naturally arise from our results.

\subsection{Questions about Dehn twists along 3-spheres}

Theorem~\ref{thm: elliptic surface sum}, in particular, shows that
the Dehn twist along the neck of
\[
E(8k+2)\#E(8l+2)
\]
is smoothly nontrivial for all $k,l\geq0$.
It is natural to ask whether the same conclusion holds for arbitrary simply-connected spin elliptic surfaces.
A particularly natural case is the following:

\begin{ques}
Let $k,l\geq1$ be integers.
Is the Dehn twist along the neck of $E(2k)\#E(2l)$ not smoothly isotopic to the identity?
\end{ques}

The restriction to spin elliptic surfaces (i.e., to $E(n)$ with even $n$) is essential.
Indeed, if $n$ is odd, the boundary Dehn twist on $E(n)^\circ$ is smoothly isotopic to the identity by~\cite[Corollary~A.5]{OrsonPowell2022}, since $E(n)$ is simply-connected and non-spin.
Consequently, the Dehn twist along the neck of $E(n)\#E(m)$ is smoothly isotopic to the identity whenever at least one of $n$ and $m$ is odd.

Another natural question concerns stabilization.
Lin~\cite{Lin20Dehn} proved that the Dehn twist along the neck of $K3\#K3$ remains smoothly nontrivial after a single stabilization by $S^2\times S^2$.
It is therefore natural to ask whether the nontriviality established in Theorem~\ref{thm: elliptic surface sum} persists after a single stabilization.
A particularly natural case is the following:

\begin{ques}
Do there exist integers $k,l\geq0$ with $(k,l)\neq(0,0)$ such that the diffeomorphism $\delta\#\id_{S^2\times S^2}$ of
\[
E(8k+2)\#E(8l+2)\#(S^2\times S^2)
\]
is not smoothly isotopic to the identity?
Here, $\delta$ denotes the Dehn twist along the neck of $E(8k+2)\#E(8l+2)$.
\end{ques}

A major challenge in extending Lin's argument to our setting is that the families Bauer--Furuta invariants associated with the spin$^c$ structures considered in this paper, which need not arise from spin structures, do not generally admit the $\Pin(2)$-symmetry that played a crucial role in~\cite{Lin20Dehn}.

\subsection{Questions about Seifert fibered Dehn twists}

By \cite[Proposition~2.25]{KLMME}, the boundary Dehn twist on the standard nucleus $N(n)=N(n)_{1,1}$ is smoothly isotopic to the identity for every $n\geq2$.
On the other hand, Corollary~\ref{cor: nuclei} shows that, when $n$ is divisible by 4, the boundary Dehn twist on $N(n)_{p,q}$ is smoothly nontrivial for odd coprime $p,q$ with $(p,q)\neq(1,1)$.
One may ask whether the following generalization holds:

\begin{ques}
Let $n\geq2$ be an integer, and let $p,q$ be relatively prime positive integers.
Is the boundary Dehn twist on $N(n)_{p,q}$ smoothly isotopic to the identity relative to the boundary if and only if $(p,q)=(1,1)$?
\end{ques}

%\begin{ques}
%Let $n\geq2$ be an even integer, and let $p,q$ be relatively prime positive odd integers.
%Is the boundary Dehn twist on $N(n)_{p,q}$ smoothly isotopic to the identity relative to the boundary if and only if $(p,q)=(1,1)$?
%\end{ques}

Another natural question concerns the orders of Dehn twists along Seifert fibered 3-manifolds.
Unlike Dehn twists along 3-spheres, which have order at most two, Dehn twists along Seifert fibered 3-manifolds can have infinite order, as shown in~\cite{KangParkTaniguchi,KLMME}.
This motivates the following question:

\begin{ques}
Let $k\geq1$, and let $p,q$ be relatively prime positive odd integers with $(p,q)\neq(1,1)$.
Does the boundary Dehn twist on $N=N(4k)_{p,q}$ have infinite order in $\pi_0(\Diff_\del(N))$?
\end{ques}

Note that our proof of Corollary~\ref{cor: nuclei} does not rule out the possibility that some even power of the boundary Dehn twist on $N(4k)_{p,q}$ is smoothly isotopic to the identity.
Indeed, the proof establishes the nontriviality of the boundary Dehn twist by extending it to $E(4k)_{p,q}$ and applying the homomorphism
\[
\FBF_{\fraks}^{\divi}: \pi_0(\TDiff(E(4k)_{p,q})) \to \pi_2^s\cong\Z/2.
\]
Since the target is $\Z/2$, this homomorphism vanishes on all even powers of any diffeomorphism.

\bibliographystyle{alpha}
\bibliography{mainref}

\end{document}